\documentclass[a4paper,11pt]{article}

\usepackage{mathbbold}
\usepackage{amsfonts}
\usepackage{epsf,latexsym,amsfonts,amsbsy,mathrsfs}
\usepackage{cite}
\usepackage{amsmath, amssymb, amsthm,amscd,amsxtra}
\usepackage[]{graphicx,subfigure}
\usepackage{float}
\usepackage{multirow}
\usepackage{bbm}
\usepackage{stmaryrd}
\usepackage{epsfig}

\newtheorem{theorem}{Theorem}[section]

\newtheorem{algorithm}{Algorithm}[section]

\theoremstyle{definition}
\newtheorem{definition}[theorem]{Definition}

\theoremstyle{remark}

\numberwithin{equation}{section}

\newcommand{\ii}{\rm i}

\newcommand{\pd}[1]{\left\langle #1\right\rangle}

\newcommand{\nn}{\nonumber}

\newcommand{\Bx}{{\boldsymbol{x}}}

\newcommand{\BA}{{\boldsymbol{A}}}

\newcommand{\BD}{{\boldsymbol{D}}}

\newcommand{\BI}{{\boldsymbol{I}}}

\newcommand{\BL}{{\boldsymbol{L}}}
\newcommand{\BM}{{\boldsymbol{M}}}

\newcommand{\BR}{{\boldsymbol{R}}}
\newcommand{\BS}{{\boldsymbol{S}}}

\newcommand{\BU}{{\boldsymbol{U}}}

\newcommand{\Ce}{{\mathcal E}}

\newcommand{\Cg}{{\mathcal G}}

\newcommand{\Ct}{{\mathcal T}}

\newtheorem{rem}{Remark}[section]
\newtheorem{exm}{Example}[section]

\title{Multilevel Preconditioner with Stable Coarse Grid Corrections for the Helmholtz Equation}

\author{Huangxin Chen\thanks{School of Mathematical Sciences, Xiamen University, Fujian, 361005, P.R. China ({\tt chx@xmu.edu.cn}).},
Haijun Wu\thanks{Department of Mathematics, Nanjing University,
Jiangsu, 210093, P.R. China ({\tt hjw@nju.edu.cn}).}, and Xuejun
Xu\thanks{Institute of Computational Mathematics and
Scientific/Engineering Computing, Academy of Mathematics and Systems
Science, Chinese Academy of Sciences, P.O. Box 2719, Beijing,
100190, P.R. China ({\tt xxj@lsec.cc.ac.cn}).} }

\date{}

\begin{document}
\maketitle

\begin{abstract}
In this paper we consider a class of robust multilevel precontioners
for the Helmholtz equation with high wave number.  The key idea in
this work is to use the continuous interior penalty finite element
methods (CIP-FEM) studied in \cite{Wu12,Wu12-hp} to construct the
stable coarse grid correction problems. The multilevel methods,
based on GMRES smoothing on coarse grids, are then served as a
preconditioner in the outer GMRES iteration. In the one dimensional
case, convergence property of the modified multilevel methods is
analyzed by the local Fourier analysis. From our numerical results,
we find that the proposed methods are efficient for a reasonable
range of frequencies. The performance of the algorithms depends
relatively mildly on wave number. In particular, only one GMRES
smoothing step may guarantee the optimal convergence of our
multilevel algorithm, which remedies the shortcoming of the
multilevel algorithm in \cite{EEO01}.

\end{abstract}

\noindent {\bf Keywords.} Multilevel method, Helmholtz equation,
high wave number, continuous penalty finite element method, GMRES
method, local Fourier analysis

\section{Introduction}

The efficient and accurate numerical approximation of high frequency
wave propagation is of fundamental importance in many applications
such as acoustic, electromagnetic, elasticity and geophysical
surveys. When the problem is linear and time-harmonic, it can be
typically modeled by Helmholtz equation. The interest of this paper is to
consider the Helmholtz equation with Robin boundary condition which is
the first order approximation of the radiation condition:
\begin{align}
\label{Hem}
-\triangle u-\kappa^2u&= f  \qquad {\rm in }\ \Omega,\\
\label{Robin} \frac{\partial u}{\partial {n}}+{\bf i}\kappa u&=  g
\qquad {\rm on }\ \partial \Omega,
\end{align}
where $\Omega \subset \mathbb R^d,d=2,3$ is a polygonal/polyhedral domain, $\kappa>0$
is known as the wave number, ${\bf i}=\sqrt{-1} $ denotes the
imaginary unit, and ${n}$ denotes the unit outward normal to
$\partial \Omega$.

For high wave number $\kappa$, the linear system from the
discretization of Helmholtz equation is usually strongly indefinite,
causing most of iterative methods to converge slowly or diverge. In
recent years, there have been many advances in the development of
iterative methods and preconditioners for the solution of the
Helmholtz equation (cf. \cite{Erlangga08,Ernst11}).

Due to high efficiency of multilevel methods for positive definite
problems, more and more attentions have been received to develop
robust multilevel methods for Helmholtz equation. However, a direct
application of multilevel methods with standard smoothing and coarse
grid correction are ineffective. Some strategies have been proposed
to remedy the problem (cf. \cite{Brandt97,EEO01,Kim,Livshits}). For
instance, Elman, Ernst and O'Leary \cite{EEO01} proposed GMRES
smoothing together with flexible GMRES acceleration. But in order to
achieve convergence, a relatively large number of GMRES smoothing
steps are needed on the intermediate grids. Another approach, the
so-called wave-ray multigrid methods \cite{Brandt97,Livshits}, was
proposed by Brandt and Livshits through defining a meaningful coarse
problem by augmenting the standard V-cycle with ray grids and using
coarse grid basis derived from plane waves. This method performs
well with increasing wave number, but it does not easily generalize
to unstructured grids and complicated Helmholtz problems.
Alternatively, instead of applying multigrid iterations directly to
the Helmholtz equation, a class of shifted Laplacian preconditioners
\cite{Erlangga04,Erlangga06} has recently attracted a lot of
attention, which precondition the Helmholtz equation with a complex
shifted operator and is shown to be an efficient Krylov method
preconditioner. We would like to mention that Engquist and Ying
\cite{Engquist1,Engquist2} recently proposed two new types of
sweeping preconditioners for central difference scheme of the
Helmholtz equation based on an approximate $LDL^t$ factorization by
sweeping the domain layer by layer starting from an absorbing layer.
Similar to the wave-ray multigrid methods, the new preconditioners
have a nearly linear computational cost and the number of outer
iterations is essentially independent of the number of unknowns and
the wave number when combined with the GMRES solver.

Our objective is to develop robust multilevel methods for the
Helmholtz problem (\ref{Hem})-(\ref{Robin}). Although the pollution
error is inherent in the standard finite element methods (FEM) to
solve Helmholtz equation, the FEM can still be used on fine grid
whose mesh size is sufficiently small to reduce the pollution error.
In a recent work \cite{Wu12-hp}, the pre-asymptotic analyses of both
the FEM and  the continuous interior penalty finite element methods
(CIP-FEM) are given. In particular, the well-posedness of standard
FEM has been proved under the condition of $\frac{\kappa h}{p}\leq
C_0 (\frac{p}{\kappa})^{\frac{1}{p+1}}$, where $h$ is mesh size, $p$
is the polynomial order of approximation space, and $C_0$ is a
constant independent of $\kappa,h,p$. Thus, without this condition,
the well-posedness of standard FEM on coarse grids in the multilevel
method can not be guaranteed. Moreover, oscillations on the scale of
the wavelength can not be resolved well by standard FEM on the
coarse grids. By contrast, the CIP-FEM has been proved in
\cite{Wu12-hp} that the associated discrete problem is always
well-posed without any mesh condition. Intrinsically, the main
technique in the stable CIP-FEM is to add a complex shift in the
bilinear form, which is similar to the idea of shifted Laplace
preconditioner approach. Comparing to adding a shift to the original
problem in shifted Laplace operator, the well-posed CIP-FEM is
consistent with the original equation and only changes the discrete
bilinear form. Based on these observations, the new approach
proposed in this work is to apply the CIP-FEM to construct the
stable coarse grid correction problems. Standard Jacobi or
Gauss-Seidel smoothers become unstable on the coarse grids,
especially on the intermediate grids in multilevel methods.
Motivated by the smoothing approach presented in \cite{EEO01}, the
smoothing in this work is to use GMRES method based on CIP-FEM on
the coarse grids, and standard Jacobi or Gauss-Seidel relaxation on
the fine enough grids. From our numerical experiments, we find that
the number of GMRES smoothing steps in our algorithm can be much
smaller than that in \cite{EEO01}, even if one GMRES smoothing step
may guarantee the optimal convergence of our multilevel algorithm.

Our main tool to analyze the multilevel methods for Helmholtz
equation is the {\it Local Fourier analysis} (LFA), which has been
introduced for multigrid analysis by Achi Brandt in 1977
\cite{Brandt77}. Comprehensive surveys can be found in \cite{WJ04}
and the references therein. We mainly utilize the LFA to analyze
smoothing properties of relaxations and convergence properties of
two- and three-level methods in one dimensional case. This may
provide quantitative insights into the multilevel methods for
Helmholtz problem (\ref{Hem})-(\ref{Robin}).

The remaining part of this paper is organized as follows: In section \ref{notation}, we introduce some
notation, recall the formulation of CIP-FEM, and present the multilevel method for the linear system from
CIP-FEM approximation. Section \ref{modified-multilevel} is to present the modified multilevel method for Helmholtz problem.
Standard FEM is used on fine enough grids, on coarse grids the CIP-FEM is utilized instead. Section \ref{LFA} is devoted
to the LFA of the multilevel method for one dimensional Helmholtz problem, we focus on the smoothing analysis and
two- and three-level analysis. In the last section, we give some numerical results to illustrate the
performance of the proposed multilevel methods.

\section{Notations and Preliminaries}\label{notation}
\subsection{Formulation of CIP-FEM}
Let $\mathcal T_h$ be a conforming quasi-uniform triangulation of $\Omega$, and denote the collection of
edges/faces by $\mathcal E_h$, while the set of interior edges/faces by
$\mathcal E_h^I$ and the set of boundary edges/faces by $\mathcal E_h^B$.
For any $T \in \Ct_h$, we define $h_T := {\rm diam}(T)$.
Similarly, for $e \in \Ce_h$, define $h_e:={\rm diam}(e)$. Let
$h:=\max_{T \in \Ct_h}h_T$. Throughout this paper we use the
standard notations and definitions for Sobolev spaces (cf.
\cite{Adams}). In particular, $(\cdot,\cdot)_Q$ and $\pd{ \cdot,\cdot}_\Sigma$
for $\Sigma\subset \partial Q$ denote the $L^2$-inner product
on complex-valued $L^2(Q)$ and $L^2(\Sigma)$
spaces, respectively. Denote by $(\cdot,\cdot):=(\cdot,\cdot)_\Omega$
and $\pd{ \cdot,\cdot}:=\pd{ \cdot,\cdot}_{\partial\Omega}$.

Now we define the energy space $V:= H^1(\Omega) \cap \prod_{T \in
\Ct_h} H^2(T) $. For any $v \in V$ and an interior edge/face $e = T_1
\cap T_2$, where $T_1$ and $T_2$ are two distinct elements of
$\Ct_h$ with respective outer normals $n_1$ and $n_2$, we introduce
the jump $\Lbrack \nabla v \cdot n \Rbrack|_e = \nabla v \cdot
n_1|_{T_1} + \nabla v \cdot n_2|_{T_2} $. Define the sesquilinear
form $b_h(\cdot,\cdot)$ on $V \times V$ as follows:
\[
b_h(u,v) := (\nabla u, \nabla v) + J(u,v),\qquad  u,v \in V,
\]
where
\begin{align}
J(u,v) := \sum_{e \in \Ce_h^I} {\bf i} \gamma_eh_e
\pd{\Lbrack \nabla u \cdot n \Rbrack,\Lbrack \nabla v \cdot n
\Rbrack}_e,\label{Jump-var}
\end{align}
where ${\bf i}\gamma_e$ is a complex
number with positive  imaginary part. The terms in $J(u,v)$ are
so-called penalty terms and ${\bf
i}\gamma_e$ are  penalty parameters (cf. \cite{Wu12,Wu12-hp}).

Clearly, $J(u,v)=0$ if $u\in H^2(\Omega)$ and $v \in V$. Thus, if
$u\in H^2(\Omega)$ is the solution of (\ref{Hem})-(\ref{Robin}),
then there holds
\begin{equation}
a_h(u,v):=b_h(u,v) -\kappa^2 (u,v) +{\bf i} \kappa \langle u,v\rangle = (f,v)
+\langle g,v\rangle, \qquad  v \in V.\label{variation}
\end{equation}

We define the  CIP approximation space $V_h$ as the standard  finite element space of order $p$, i.e.,
\[
V_h := \big\{ v_h \in H^1(\Omega): v_h |_T \in \mathcal {P}_p(T),
\,T \in \Ct_h \big\},
\]
where $\mathcal {P}_p(T)$ is the space of polynomials of degree at
most $p$ on $T$. The CIP finite element
approximation is to find $u_h \in V_h$ such that
\begin{equation}
a_h(u_h,v_h)
= (f,v_h) +\langle g,v_h\rangle, \qquad v_h
\in V_h.\label{variation}
\end{equation}

It is clear that if the parameters $\gamma_e\equiv 0$, then the CIP-FEM
reduces to the standard  FEM. It has been proved
that the CIP-FEM is stable for any $\kappa,h,p>0$ \cite{Wu12,Wu12-hp}. The penalty parameters
may be tuned to reduce the pollution errors. The numerical results
in \cite{Wu12} show that using about the same total degrees of
freedom (DOFs), the CIP-FEM yields the least phase error comparing
to the standard FEM and  IPDG method \cite{Wu09}.

\subsection{Multilevel methods for CIP-FEM}

Let $\{\mathcal {T}_{l}\}_{l=0}^L$ be a shape regular family of
nested geometrically conforming simplicial triangulations of the
computational domain $\Omega$ obtained by successive quasi-uniform
refinement of an intentionally chosen coarse grid $\mathcal
{T}_{0}$. We denote by $V_l$ the CIP approximation space on $\Ct_l$.
It is easy to see that the spaces $\{V_l\}_{l=0}^L$ are nested,
i.e., $V_0 \subset V_1 \subset \cdots \subset V_L$. But the bilinear
forms $\{a_l(\cdot,\cdot)\}_{l=0}^L$ defined as in (\ref{variation})
on each $V_l$ are nonnested. Thus in this work we consider the
following multilevel methods for CIP-FEM discretizations, which has
also been used for solving the linear systems from nonconforming P1
finite element approximations (cf. \cite{VW}).

For brevity, we denote by $a_l(\cdot,\cdot)$ the bilinear form
$a_{h_l}(\cdot,\cdot)$ on $V_l$, where $h_l$ is the mesh size of
$\Ct_l$. Define projections $P^C_l$, $Q_l$ : $V_L \rightarrow V_l$
according to
$$
a_l(P^C_l v, w) = a_L(v,  w)   , \quad (Q_lv,w) = (v, w)  , \quad\forall
v\in V_L,\, w \in V_l.
$$
The existence and uniqueness  of the discrete problem (\ref{variation})
imply the well-posedness of each $P^C_l$. For $0\leq l \leq L$, we
also define $A^C_l: V_l \rightarrow V_l$ by means of
\begin{align}
(A^C_l v,w) = a_l(v,w), \quad \forall v, w\in V_l.\label{Ah}
\end{align}
Define $F_l\in V_l$ by $$(F_l,v) =
(f,v) + \langle g,v\rangle\quad \forall v\in V_l.$$
Then the CIP-FEM on level $l$ (cf. \eqref{variation}) can be rewritten as:
\begin{equation}\label{ecipfeml}
A^C_lu^C_l=F_l.
\end{equation}

For the smoothing strategy, in fact, when $\kappa h_l/p $ is small
enough, either of weighted Jacobi and Gauss-Seidel relaxation can be
applied. Otherwise, GMRES relaxation can be used as a smoother and
this will be explained in the following sections. To be precise, we
describe the smoothing strategy as follows:
\begin{definition}\label{dRlC} Let $R_0^C=(A_0^C)^{-1}$. Given $\alpha>0$ and $0<l\le
L$, let $S_L=\{ l : \ \kappa h_l / p < \alpha, \ 1 \le l\le L\}$ and
$G_L = \{l : 1 \le l \le L\} \setminus S_L$. If $l \in S_L$, let
$R_l^C$ be the weighted Jacobi relaxation $R^{J}_{l,C}$ or
Gauss-Seidel relaxation $R^{GS}_{l,C}$ based on $A^C_l$. Otherwise,
when $l \in G_L$ we use the GMRES relaxation based on $A^C_l$.
\end{definition}
We remark that we choose $\alpha=0.5$ in this paper. This choice is
motivated by the efficient relaxation of Jacobi and Gauss-Seidel
smoothers, and it ensures that the amplification factor for the
Jacobi and Gauss-Seidel smoothers will not become too large (cf.
section 4.3 in the following and section 2.1 in \cite{EEO01}). When
$l \in S_L$, we also define the operator $(R_l^C)^t$  by
\begin{align}
\big((R_l^C)^t\bar w,\bar v\big)=\big(R_l^Cv,w\big),\quad\forall
v,w\in V_l. \label{sm-op}
\end{align}

We can now state the multilevel method for solving the CIP-FEM
discretization system on level $L$ which is non-recursive version of
multilevel method.

 \smallskip

\begin{algorithm}\label{alg1}
Given an arbitrarily chosen initial iterate $u^0\in V_L$, we seek
$u^n\in V_L$ as follows:

\smallskip

Let $v_0 = u^{n-1}$.

\begin{itemize}
\item[$1)$] For $l=0,1,\cdots,L$. When $l=0$ or $l \in S_L$,
\[
v_{l+1} = v_l + \mu_l R_l^CQ_l(F_L - A^C_L v_l).
\]
Otherwise, perform GMRES relaxation for the correction problem
$A^C_l w_l = Q_l(F_L - A^C_L v_l)$ and set $v_{l+1} = v_l + \mu_l
w_l$. Here $\mu_{l} >0$ is a scaling parameter to weaken the
influence of the error during prolongations. In this paper, we
always set $\mu_{l}\equiv 0.5$.

\item[$2)$] For $l=L,\cdots,1,0$. When $l=0$ or $l \in S_L$,
\begin{align*}
v_{2L+2-l} = v_{2L+1-l} + \mu_l (R_l^C)^tQ_l(F_L - A^C_L
v_{2L+1-l}).
\end{align*}
Otherwise, perform GMRES relaxation based on $v_{2L+1-l}$ as in step
1 to get $v_{2L+2-l}$.

\item[$3)$] Set $u^n = v_{2L+2} $.
\end{itemize}
\end{algorithm}

Here we note that the GMRES relaxation is not linear in the starting
value, and the error operator of the above algorithm can not be
written directly in a product form which holds only for the
particular case $G_L = \emptyset$. When this particular case is
considered reasonably, we set
\begin{align*}
T_l := \mu_{l} R_l^C A^C_l P^C_l\text{ and } T_l^* := \mu_{l}
(R_l^C)^t A^C_l P^C_l ,   \quad l=0,1,\ldots,L, \ G_L = \emptyset,
\end{align*}
Then the error operator of Algorithm \ref{alg1} for the case
$G_L=\emptyset$ can be derived as
\begin{align}
E_M^*E_M, \text{ where } E_M := (I-T_L)\cdots(I-T_1)(I-T_0) , E_M^* := (I-T_0^*)(I-T_1^*)\cdots(I-T_L^*), \label{EM-ML}
\end{align}
where $I$ is the identity operator in $V_L$.


\section{Modified multilevel methods for standard FEM}\label{modified-multilevel}

In general, in order to reduce the pollution error, 6-10 grid points
per wavelength are typically chosen to yield reasonable accuracy. In
a unit square domain $[0,1]^2$, it is well know that $\kappa h =
2\pi/n_w$, where $n_w$ is the number of points per wavelength (cf.
\cite{Ihlenburg_book}), implies the deterioration of $\kappa h$ for
increasing grid points per wavelength. In theory, the well-posedness
of discrete solution by CIP-FEM holds for any $\kappa,h$ and $p$,
but it is not guaranteed  for discrete solution by standard FEM on
coarser grids with $\kappa h/p \geq C$ particularly
\cite{Wu12,Wu12-hp}, where the constant $C$ is independent of
$\kappa,h$ and $p$. Therefore, in order to establish stable coarse
grid correction problems in multilevel methods for Helmholtz
equation, the CIP-FEM will be applied on the coarse grids. Besides,
for standard FEM, eigenvalues of discrete system close to the origin
may undergo a sign change after discretization on a coarser grid. If
a sign change occurs, the coarse grid correction does not give a
convergence acceleration to the finer grid problem but gives a
severe convergence degradation instead. This is analyzed in
\cite{EEO01} and a remedy combining GMRES method as a smoother on
coarse grids is proposed. This idea has been applied in Algorithm
\ref{alg1}, and we will also utilize this strategy in the following
modified multilevel methods to get more efficient smoother for
indefinite discrete systems.

\smallskip

For brevity, let $A^F_l$ denote the linear operator $A^C_l$ with the
parameters $\gamma_e\equiv0$, i.e., the linear operator for standard
FEM, and energy operator $P^F_l = P^C_l |_{\gamma_e \equiv 0}$. The
discrete problem on level $l$ is $A^C_lu^C_l=F_l$ for CIP-FEM (cf.
\eqref{ecipfeml}) or $A^F_lu^F_l=F_l$ for standard FEM. When the
mesh size is sufficiently small to reduce the pollution error and
satisfy the accuracy requirement, both CIP-FEM and standard FEM can
be utilized. However, the nonzero elements of linear system from
standard FEM may be less than that from CIP-FEM, and when the grid
is fine enough, the pollution error by the standard FEM is also
small. Thus, we may still apply the standard FEM on the fine grids.

\smallskip

Similar to Definition \ref{dRlC}, we use the smoothing strategy on
each $V_l$ as follows.
\begin{definition}\label{dRl} Let $R_0=(A_0^C)^{-1}$. Given $\alpha>0$ and $0<l\le L$.  If
$l \in S_L$, let $R_l$ be the weighted Jacobi relaxation
$R^{J}_{l,F}$ or Gauss-Seidel relaxation $R^{GS}_{l,F}$ based on
$A^F_l$. Otherwise, when $l \in G_L$ we use the GMRES relaxation
based on $A^C_l$.
\end{definition}
When $l \in S_L$, the operator $R_l^t$ is defined similarly as in
(\ref{sm-op}). Now we state the modified multilevel method for the
discrete system from standard FEM for the Helmholtz problem
(\ref{Hem})-(\ref{Robin}) as follows.

 \smallskip

 \begin{algorithm}\label{modi-mu}
Given an arbitrarily chosen initial iterate $u^0\in V_L$ and
integers $m_1, m_2\ge 0$, we seek $u^n\in V_L$ as follows:

\smallskip

Let $v_0 = u^{n-1}$.

\begin{itemize}
\item[$1)$] For $l=0,1,\cdots,L$. Given $\mu_{l} >0$ is chosen
as in Algorithm \ref{alg1}, when $l=0$ or $l \in S_L$,
\[
v_{l+1} = v_l + \mu_l (R_l)^{m}Q_l(F_L - A^F_L v_l),
\]
where $m=m_1$ if $l>0$, $m=1$ if $l=0$. Otherwise, perform $m_1$
steps of GMRES relaxation for the correction problem $A^C_l w_l =
Q_l(F_L - A^F_L v_l)$ and set $v_{l+1} = v_l + \mu_l w_l$.

\item[$2)$] For $l=L,\cdots,1,0$. When $l=0$ or $l \in S_L$,
\begin{align*}
v_{2L+2-l} = v_{2L+1-l} + \mu_l (R_l^t)^mQ_l(F_L - A^F_L
v_{2L+1-l}),
\end{align*}
where $m=m_1$ if $l>0$, $m=1$ if $l=0$. Otherwise, perform $m_2$
steps of GMRES relaxation for the correction problem $A^C_{l} w_{l}
= Q_{l}(F_L - A^F_L v_{2L+1-l})$ and set $v_{2L+2-l} = v_{2L+1-l} +
\mu_{l} w_{l}$.

\item[$3)$] Set $u^n = v_{2L+2} $.
\end{itemize}
\end{algorithm}

\begin{rem}
In \cite{EEO01}, in order to prevent unnecessary damping of
smoothing modes which should be handled by the coarse grid
correction, the postsmoothing is always favored over presmoothing
(cf. section 3.2 in \cite{EEO01}). This is also true in our work.
However, due to the utilization of stable coarse grid correction
method, the above algorithm will converge even when the smoothing is
chosen as one step in both post- and presmoothing.
\end{rem}

\begin{rem}

Comparing to Algorithm \ref{alg1}, the CIP-FEM is only applied in
the coarse grid correction when $l \in G_L$ in the above algorithm.
From the first numerical example in this paper, we can see that the
convergence property of this two algorithms is similar. Actually,
this phenomena can also be observed in the following LFA. In order
to reduce computations, one may prefer Algorithm \ref{modi-mu} in
practice.

\end{rem}

\section{The one dimensional local Fourier analysis}\label{LFA}

In this section, we aim to analyze different approaches based on
Algorithm \ref{alg1} for the discrete system from one dimensional
Helmholtz equation, where standard FEM is utilized on the finest
grid. We focus on the analysis for the discretization from linear
continuous interior penalty finite element method (CIP-P1) by the
important tool LFA in multigrid analysis. Here we mainly focus on
the analysis of two-level methods. The LFA of three-level methods is
also mentioned. The analysis will imply the efficiency of Algorithm
\ref{modi-mu}. The following presentation is related to the
notation and philosophy from \cite{WJ04}.

\subsection{Basic tools in one dimensional local Fourier analysis}

The LFA is based on certain idealized assumptions and
simplifications: the boundary conditions are neglected and the
problem is considered on regular indefinite grids $\Cg_h = \{ x:
x=x_j=jh, j \in \mathbb{Z} \}$. Although the Robin boundary condition (\ref{Robin}) and other
absorbing boundary conditions are often applied in realistic Helmholtz problem, the neglect of
boundary condition does usually not affect the validity of LFA (cf. \cite{WJ04}). The LFA
in this section can be considered as simplification for the analysis of multilevel method
for Helmholtz problem (\ref{Hem}-\ref{Robin}) in one dimension.

Let $\BL_h$ be a discrete operator
with a stencil representation $\BL_h = [l_j]_h, j \in \mathbb{Z}$.
For any $u_h$ defined on $\Cg_h$ and a fixed point $x\in \Cg_h$,
$\BL_hu_h$ reads in stencil notation as
\begin{align}
\BL_h u_h(x) = \sum_{j \in J} l_j u_h(x+jh), \nn
\end{align}
where $J\subset \mathbb{Z}$ is a certain finite index defining the
so-called stencil. The basic quantities in the LFA are the Fourier modes $
\varphi_h(\theta,x) = e^{{\bf i} \theta x/h}$ with $x \in \Cg_h$ and
Fourier frequency $\theta \in \mathbb{R}$. In fact, the frequency
$\theta$ can be restricted to the interval $(-\pi,\pi]$ as a fact
that $\varphi_h(\theta + 2\pi,x) = \varphi_h(\theta,x)$. It is easy
to see that the Fourier modes are all the formal eigenfunctions of
$\BL_h$:
\begin{align}
\BL_h \varphi_h(\theta,x) = \widetilde{\BL}_h(\theta)
\varphi_h(\theta,x), \qquad x\in \Cg_h , \theta \in (-\pi,\pi],
\label{Lh-stencil}
\end{align}
where the eigenvalues of $\BL_h$ can be presented as $
\widetilde{\BL}_h(\theta) = \sum_{j\in J} l_j e^{{\bf i} \theta j}$,
which is called the Fourier symbol of the operator $\BL_h$. Given a
so-called low frequency $\theta^0\in \Theta_{\rm
low}=(-\pi/2,\pi/2]$, its complementary frequency $\theta^1$ is
defined as
\begin{align}
\theta^1 = \theta^0-{\rm sign}(\theta^0)\pi. \label{theta1}
\end{align}
Interpreting the Fourier modes as coarse grid functions yields
\[
\varphi_h(\theta^0,x) = \varphi_{2h}(2\theta^0,x) =
\varphi_{2h}(2\theta^1,x) = \varphi_h(\theta^1,x), \qquad \theta^0
\in \Theta_{\rm low}, \ x\in \Cg_{2h}.
\]
The Fourier modes $\varphi_h(\theta^0,x) $ and
$\varphi_h(\theta^1,x) $ are called $2h$-harmonics. These Fourier
modes coincide on the coarse grid with mesh size $H=2h$, and they
can be represented by a single coarse grid mode
$\varphi_{2h}(2\theta^0,x)$. Hence, each low frequency mode is
associated with a high frequency mode. For a given $\theta^0 \in
\Theta_{\rm low}$, define the two dimensional subspace of
$2h$-harmonics by
\begin{align}
E_{2h}^{\theta^0} := {\rm span} \{
\varphi_h(\theta^0,x),\varphi_h(\theta^1,x)  \}, \label{Etheta}
\end{align}
where $\theta^1$ is defined as in (\ref{theta1}). A crucial
observation is that the space $E_{2h}^{\theta^0}$ is invariant under
both smoothing operators and correction schemes for general cases by
two-level method. The invariance property holds for many well-known
smoothing methods (cf. \cite{WJ04}), such as Jacobi relaxation,
lexicographical Gauss-Seidel relaxation, et al.

Let $\BM_h$ be a discrete two-level operator. In the following, we
will show that a block-diagonal representation for $\BM_h$ consists
of $2\times 2$ blocks $\widetilde{\BM}_h(\theta)$ (cf. \cite{WJ04}),
which denotes the representation of $\BM_h$ on $E_{2h}^{\theta^0}$.
Then the convergence factor of $\BM_h$ by the LFA is defined as
follows:
\[
\rho(\BM_h) = {\rm sup} \{ \rho(\widetilde{\BM}_h(\theta)): \ \theta
\in \Theta_{\rm low} \},
\]
where $ \rho(\widetilde{\BM}_h(\theta))$ is the spectral radius of
the matrix $\widetilde{\BM}_h(\theta)$. We can refer to \cite{WJ04}
for generalizations to $k$-level analysis.

\subsection{One dimensional Fourier symbols}

Now we give the Fourier symbols of different operators in multilevel
method for CIP-P1 discretization of one dimensional Helmholtz equation
(\ref{Hem}). We always assume $\gamma_e \equiv\gamma$ for some
constant $\gamma$ in (\ref{Jump-var}). Denote by $t=\kappa h$, $R = -1-{\bf
i}4\gamma - t^2/6$, $S=1+{\bf i}3\gamma - t^2/3$. Since the boundary
condition is neglected in the LFA, the stencil presentation of
discretization operator $\BA^C_h$ from (\ref{Ah}) can be derived as (cf. \cite{Wu12-1d})
\begin{align}
\BA^C_h = \frac{1}{h}[ {\bf i}\gamma \ \ R  \ \ 2S \ \ R \ \ {\bf
i}\gamma ]_h.\nn
\end{align}
Combining the above expression and (\ref{Lh-stencil}) yields the
Fourier symbol of $\BA^C_h$ as
\begin{align}
\widetilde{\BA}^C_h(\theta) = \frac{1}{h}({\bf i}2\gamma
\cos{2\theta} + 2R \cos{\theta} +2S ). \label{Ah-repre}
\end{align}
Obviously, the Fourier symbol of standard FEM with piecewise P1
approximation (FEM-P1) is $\widetilde{\BA}^F_h(\theta)
=\widetilde{\BA}^C_h(\theta)|_{\gamma=0} $.

For simplicity, we use standard weighted Jacobi ($\omega$-JAC) and
lexicographical Gauss-Seidel (GS-LEX) relaxations as the smoothers
in the LFA. It is easy to derive the weighted Jacobi relaxation
matrix as $\BS^J_h = \BI_h - \omega \BD^{-1}_h \BA^C_h$, where
$\BI_h$ is indentity matrix, $\BD_h$ consists of the diagonal of
$\BA^C_h$ and $\omega$ is a weighted parameter. Due to the fact that
$\BD_h = \frac{2S}{h}\BI_h$, one can easily deduce the Fourier
symbol of weighted Jacobi relaxation as follows:
\begin{align}
\widetilde{\BS}^J_h (\theta) = 1-\frac{\omega}{S}({\bf i}\gamma
\cos{2\theta} + R \cos{\theta} +S ). \label{SJh-repre}
\end{align}
The GS-LEX relaxation matrix is $\BS^{GS}_h = (\BD_h -
\BL_h)^{-1}\BU_h$, where $-\BL_h$ is the strictly lower triangular
part of $\BA^C_h$ and $-\BU_h$ is the strictly upper triangular part
of $\BA^C_h$. The Fourier symbol of $\BS^{GS}_h$ can also be
directly derived that
\begin{align}
\widetilde{\BS}^{GS}_h (\theta) = - \frac{ R e^{{\bf i}\theta}+ {\bf
i}\gamma  e^{{\bf i}2\theta}}{ R e^{{\bf -i}\theta}+{\bf i}\gamma
e^{{\bf -i}2\theta}+2S}.\label{GSJh-repre}
\end{align}

Note that for the restriction matrix $\BI^{2h}_h=[r_j]^{2h}_h$ and
$x\in \Cg_{2h}$, there holds
\[
(\BI^{2h}_h\varphi_h(\theta^\alpha,\cdot))(x) = \sum_{j\in J}r_j
e^{{\bf i} j \theta^\alpha}\varphi_h(\theta^\alpha,x) = \sum_{j\in
J}r_j e^{{\bf i} j \theta^\alpha}\varphi_{2h}(2\theta^0,x),\qquad
\alpha = 0,1.
\]
By an analogous stencil argument, the stencil presentation of full
weighting restriction matrix can be derived to be $\BI^{2h}_h =
[1/4, 1/2, 1/4]_h^{2h}$. Then one can obtain the Fourier symbol of
$\BI^{2h}_h$ is
\[
\widetilde{\BI}^{2h}_h(\theta)  = \frac{1}{2} (1+\cos{\theta}).
\]
For the linear prolongation matrix $\BI^{h}_{2h}$, one can also
derive its Fourier symbol as follows (cf. \cite{WJ04}):
\[
\widetilde{\BI}^{h}_{2h}(\theta)  = \frac{1}{2} (1+\cos{\theta}).
\]

\subsection{Smoothing analysis}
Weighted Jacobi and lexicographical Gauss-Seidel relaxations are the
general smoothing operators. It is well-known that such two
smoothers are unstable especially for linear systems from indefinite
Helmholtz equation by standard FEM approximations. This is caused by
negative eigenvalues of the associated linear system and divergence
occurs under such two smoothers. To make up the problem, an
improvement is introduced by the shifted Laplacian preconditioner
\cite{Erlangga04,Erlangga06,Erlangga08}, which preconditions the
Helmholtz equation with a complex shifted operator
\begin{align}
-\Delta  - (1+{\bf i} \beta)k^2 , \label{shift-lap}
\end{align}
where $\beta$ are free parameter. The use of a shift is important to
guarantee multilevel convergence, whereas the multilevel method for
the shifted operator only converges for a sufficiently large shift.
But this contradicts the fact that the outer Krylov acceleration
prefers a small shift. Based on the idea of adding a shift to the
original problem, we consider the stable CIP-FEM which is consistent
with original equation and adds a complex shift in the bilinear
form.

Recalling that every two dimensional subspace of $2h$-harmonics
$E^{\theta^0}_{2h}$ with $\theta^0 \in \Theta_{\rm low}$ is left
invariant under the $\omega$-JAC and GS-LEX relaxations, then the
Fourier representation of smoother $\BS_h=\BS^J_h$ or $\BS^{GS}_h$
with respect to $E^{\theta^0}_{2h}$ can be written as
\begin{equation}
\left[
  \begin{array}{cc}
    \widetilde{\BS}_h(\theta^0)  & 0  \\
    0 & \widetilde{\BS}_h(\theta^1)\\
  \end{array}
\right],
\end{equation}
where $\widetilde{\BS}_h(\theta)$ is the smoother symbol derived in
(\ref{SJh-repre}) and (\ref{GSJh-repre}). The spectral radius of the
smoother operator can be easily calculated since the above matrix is
diagonal.

\begin{figure}[htbp]
\centering
    \includegraphics[width=2.45in]{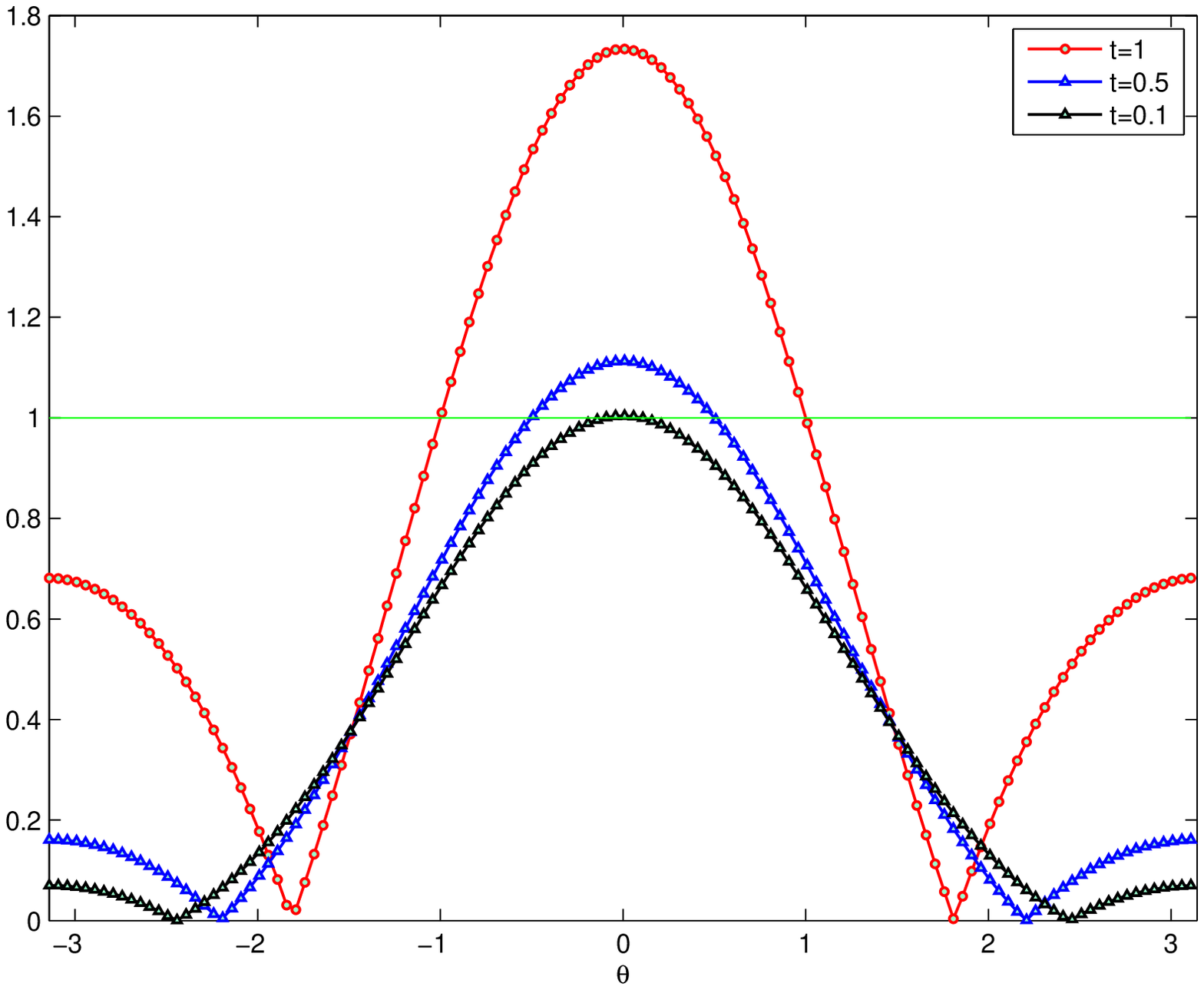}
    \includegraphics[width=2.45in]{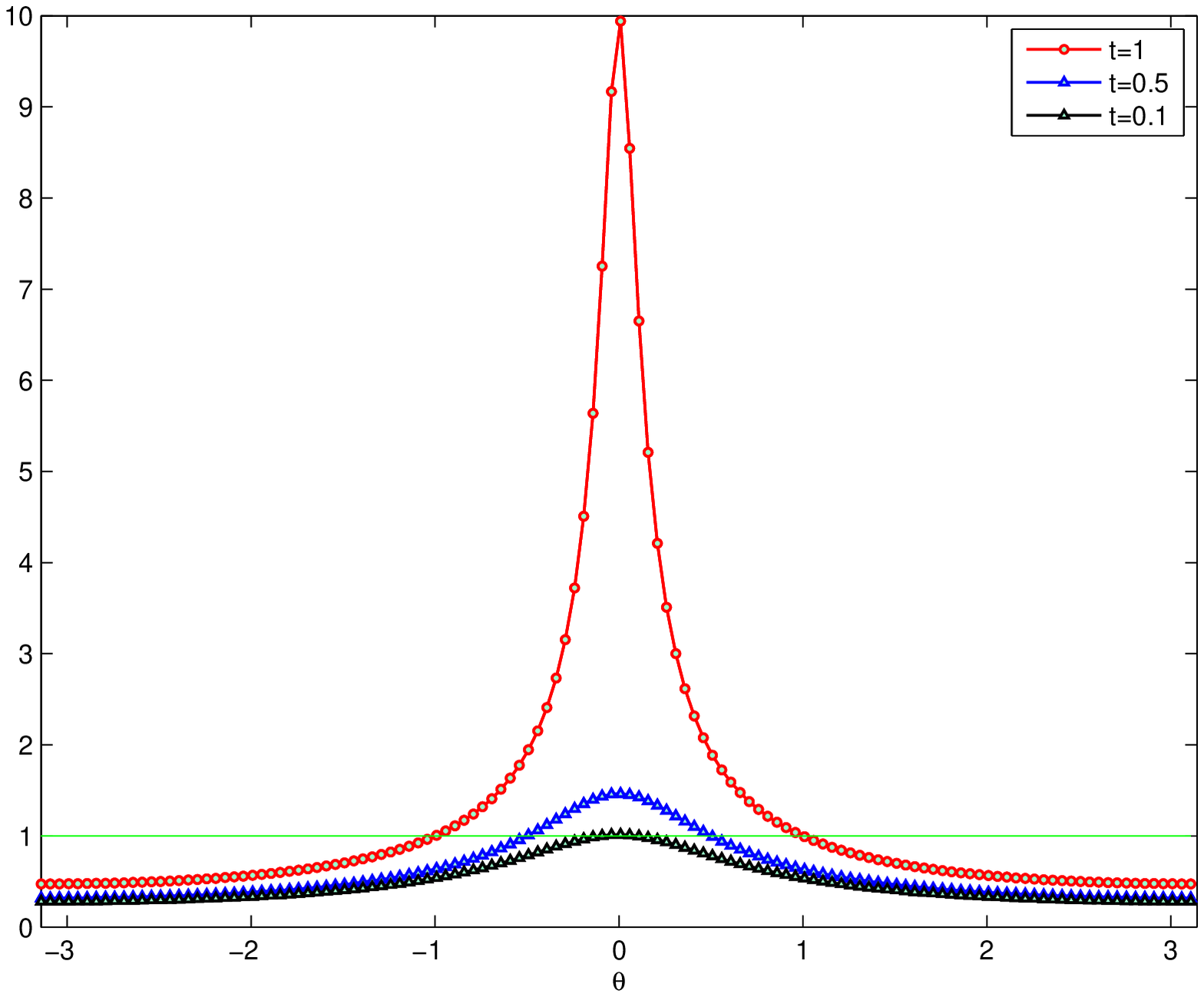}
    \caption{\small $|\widetilde{\BS}^J_h(\theta)|$ with $\omega=0.6$ (left) and $|\widetilde{\BS}^{GS}_h(\theta)|$ (right) over $(-\pi,\pi]$ for $t=0.1,0.5,1$.  }\label{fig-jacobi-1}
\end{figure}

For simplicity, we concern on the weighted Jacobi relaxation. The
frequency $\theta$ maximizing $|\widetilde{\BS}^J_h(\theta)|$ over
$(-\pi,\pi]$ can be calculated by its first and second derivatives,
which reveal $\theta=0$ or $\theta=\pi$ maximizing
$|\widetilde{\BS}^J_h(\theta)|$. Hence, the spectral radius of
$\BS^J_h$ can be deduced as follows: $ \rho(\BS^J_h) =\max\{
|\widetilde{\BS}^J_h(0)|, |\widetilde{\BS}^J_h(\pi) \}|, $ where
\begin{align*}
|\widetilde{\BS}^J_h(0) | = |1-\omega-\frac{\omega}{S}({\bf i}\gamma
+R)|, \quad |\widetilde{\BS}^J_h(\pi) | =
|1-\omega-\frac{\omega}{S}({\bf i}\gamma - R)|.
\end{align*}
Since the parameter $\gamma$ in CIP-FEM may influence the pollution
error, it is critical to make a suitable choice. For the case
$t=\kappa h \leq 1$, it has been derived in \cite{Wu12-1d} that for
one dimensional problem there exists an optimal choice ${\bf
i}\gamma_o = \frac{6\cos{t}-6+t^2\cos{t}+2t^2}{12(1-\cos{t})^2}$
such that the pollution error vanishes when the penalty parameter is
chosen as $|{\bf i}\gamma-{\bf i}\gamma_o|\leq  \frac{C}{\kappa^2
h}$, where the constant $C$ is independent of $\kappa$ and $h$. The
left graph of Figure \ref{fig-jacobi-1} shows the Fourier symbols
$\widetilde{\BS}^J_h(\theta)$ for $\omega$-JAC smoother with
$\gamma=\gamma_o$ and $\omega=0.6$. We find that
$\widetilde{\BS}^J_h(\theta)\geq 1$ always occur at the low
frequencies, and small $t$ leads to a better relaxation. Similar
phenomenon is observed for GS-LEX smoother in the right graph of
Figure \ref{fig-jacobi-1}. Thus, for fixed wave number $\kappa$,
both $\omega$-JAC and GS-LEX relaxations can be used as smoother on
fine grid. Actually, the relaxation properties of these two
smoothers are similar when ${\bf i}\gamma$ is a complex number.

When $\omega$-JAC smoother is utilized in the standard P1 FEM
approximation, then
\begin{align}
\rho(\BS^J_h)|_{\gamma=0} =\max \big\{ \mid1 -
\frac{\omega(12-t^2)}{6-2t^2}\mid, \mid 1 +   \frac{3\omega
t^2}{6-2t^2}\mid \big\}.\label{sfem-sm}
\end{align}
Thus, for $t=\kappa h$ near $\sqrt{3}$, the error is amplified by
the smoothing. For CIP-FEM, there is a (complex) shift in the
bilinear form, therefore the amplification factor is reduced. This
may permit the use of $\omega$-JAC smoother again. In fact, in our
modified multilevel method, GMRES iteration is used on the
intermediate grids (cf. \cite{EEO01}). In contrast with $\omega$-JAC
and GS-LEX relaxations, the Fourier symbol can not be derived for
GMRES smoothing. In the following, we will give some explanations
for the performance of GMRES smoothing from the numerical approach.

\begin{figure}[htbp]
\centering
    \includegraphics[width=2.5in]{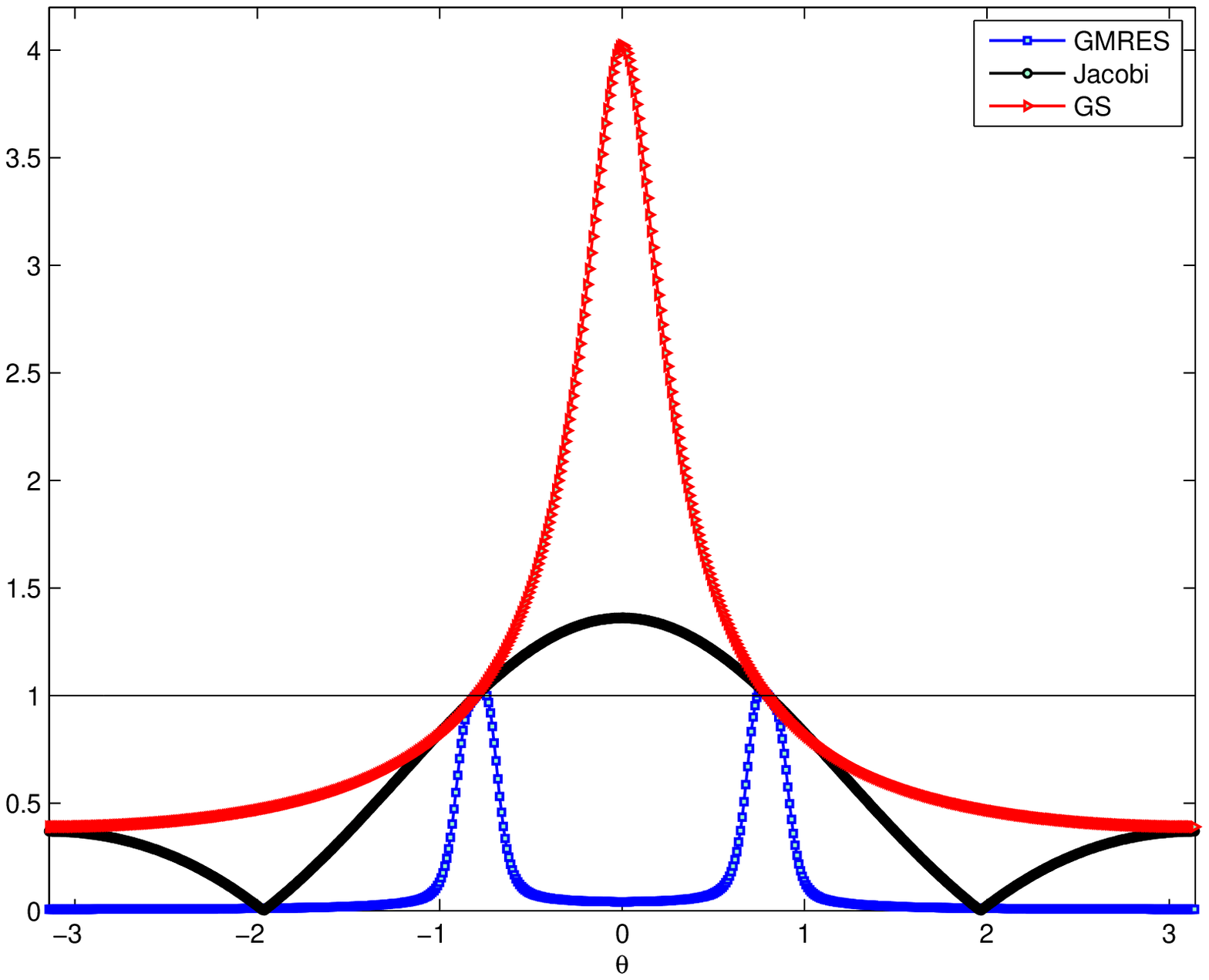}
    \caption{\small Amplification factor of GS, Jacobi and GMRES smoothing for $t=0.8$.  }\label{fig-smooth-gmres}
\end{figure}

For simplification, we consider the one dimensional Helmholtz
equation on an interval $(0,10)$ with homogeneous Dirichlet boundary
conditions. The piecewise P1 CIP-FEM discretization leads to a
linear system with $N\times N$ coefficient matrix (cf.
\cite{Wu12-1d})

\begin{eqnarray}
\BA^C_h=\frac{1}{h}\left(
\begin{array}{ccccccc}
    2S-{\bf i}\gamma &  R & {\bf i}\gamma  &  &  \\
      R & 2S & R & {\bf i}\gamma  & \\
      {\bf i}\gamma &  R & 2S & R & {\bf i}\gamma  \\
      &\ddots & \ddots & \ddots & \ddots &\ddots & \\
      & & {\bf i}\gamma &  R & 2S & R & {\bf i}\gamma  \\
      & & & {\bf i}\gamma &  R & 2S & R   \\
      & & & & {\bf i}\gamma &  R & 2S-{\bf i}\gamma \\
  \end{array} \right)_{N\times N}, \nn
\end{eqnarray}
where $N=10/h-1$. For $\kappa=200$, we apply the grid with mesh size
$h=0.004$, i.e., $t=0.8$, and choose $\gamma=\gamma_o$. Let the
vector $u_0=e^{{\bf i}\theta \Bx/h}$ be an initial choice for
smoothing, where $\Bx=[x_1,\cdots,x_{N-1}],x_k =
x_{k-1}+h,x_0=0,k=1,\cdots,N-1$, $\theta \in (-\pi,\pi]$. We assume
that $\BS_h$ is the relaxation iteration matrix and $u_1=\BS_h u_0$
is the new vector after one step of smoothing. Then for fixed
$\theta$, we obtain the amplification factor for one step of
smoothing $\rho_s(\theta)=\|u_1\| / \|u_0\|$, where $\|\cdot\|$
stands for the Euclidean norm. Figure \ref{fig-smooth-gmres} shows
the amplification factors of three different kinds of smoothing
strategies: GS, Jacobi and GMRES relaxations. We can see that for
high and low frequencies $\theta$, the amplification factor of GMRES
relaxation is always smaller than that of GS and Jacobi relaxations.
The GMRES relaxation can still be convergent even when the other two
smoothers lead to a divergent relaxation. The similar phenomenon can
also be observed for other choices of $t$ and $\gamma$. Therefore,
when standard Jacobi or GS relaxation fails as a smoother, we can
replace this with the GMRES smoothing.

From the above analysis, we can see that for the relaxation of
linear system for Helmholtz problem, the standard Jacobi or GS
relaxation can be performed on the fine grids, but they fail on the
coarse grids. In particular, the GMRES relaxation is efficient for
smoothing on the intermediate coarse grids. In the next section, the
correction scheme will be taken into account to obtain a more
realistic convergence analysis of the multilevel method.

\subsection{Two- and three-level local Fourier analysis}

For more details about convergence property of multilevel method in
Algorithm \ref{modi-mu}, we will focus on the LFA of two-level
method. The LFA of three-level method will also be mentioned. For
simplicity, we consider the case without postsmoothing. Since the
Fourier symbol can not be obtained for GMRES smoothing, we only
consider the two and three level methods with $\omega$-JAC or GS-LEX
relaxation. Then the iteration operator of Algorithm \ref{modi-mu}
in this case is the nonsymmetric version and can be derived as in
(\ref{EM-ML}). Three cases of iteration operator will be analyzed in
the following: multilevel methods by shifted Laplace approach,
multilevel methods with stable CIP-FEM coarse grid corrections, and
multilevel methods with stable CIP-FEM for fine grid smoothing and
coarse grid corrections.

The iteration operator of two-level method can be written as $
(I-T_1)(I-T_0)$, and the operator $T_l$ is chosen for the
corresponding algorithm. Since standard FEM is applied on the finest
grid, we have $T_l =\mu_l R_l A^F_l P^F_l$, where $R_l$ is smoother
determined by the algorithm. When CIP-FEM is used for smoothing on
the grids $\Ct_1$ and $\Ct_0$ respectively, the iteration matrix is
given by
\[
\BM^C_2 = \Big(\BI_1 -  \mu_1(\BI_1 - \BS^C_1)(\BA^C_1)^{-1}
\BA^F_1\Big)\Big(\BI_1 - \mu_0 \BI^1_0(\BA^C_0)^{-1} \BI^0_1
\BA^F_1\Big).
\]
Here, for $l\geq 0$, $\BS^C_l = \BI_l - \BR^C_l \BA^C_l$ is
smoothing relaxation matrix, $\BI_l$ with the same size as $\BA^F_l$
is identity matrix, $\BI^s_l(s>l)$ is prolongation matrix from level
$l$ to $s$, $\BI^s_l(s<l)$ is restriction matrix from level $l$ to
$s$, and $\BR^C_l$ stands for matrix representation of smoother
$R_l^C$ with CIP-FEM approximation. Let $\BM^{SL}_2$ (cf.
\cite{CV11}) stand for two-level iteration matrix, taking the
smoothing based on shifted Laplace operator (\ref{shift-lap}) with
standard FEM approximation. Similarly, when applying standard FEM
for smoothing on $\Ct_1$ and CIP-FEM for correction on $\Ct_0$, the
associated iteration matrix is derived as
\[
\BM^{FC}_2 = \Big(\BI_1 -  \mu_1(\BI_1 - \BS^F_1)\Big)\Big(\BI_1 - \mu_0 \BI^1_0(\BA^C_0)^{-1} \BI^0_1
\BA^F_1\Big),
\]
where $\BS^F_1  = \BI_1 - \BR^F_1 \BA^F_1$, and $\BR^F_1$ stands for
smoothing matrix representation of $R_{1}$ with standard FEM
approximation.

Since every two dimensional subspace (\ref{Etheta}) of
$2h$-harmonics $E^{\theta^0}_{2h_1}$ with $\theta^0 \in
(-\pi/2,\pi/2]$ is left invariant under $\omega$-JAC or GS-LEX
smoothing operator and correction operator, the representation of
two-level iteration matrix of $\BM^C_2$ on $E^{\theta^0}_{2h_1}$ is
given by a $2 \times 2$ matrix as follows:
\begin{eqnarray}
\widetilde{\BM}^C_2 &=& \left[ \widetilde{\BI}_1 - \mu_1\left(
\widetilde{\BI}_1 -\left[
\begin{array}{c}
    \widetilde{\BS}_1^C(\theta^0)    \\
      \widetilde{\BS}_1^C(\theta^1)\\
  \end{array} \right]_D\right) \left[  \begin{array}{cc}
    \widetilde{\BA}_1^C(\theta^0)   \\
     \widetilde{\BA}_1^C(\theta^1)\\
  \end{array}  \right]_D^{-1} \left[  \begin{array}{cc}
    \widetilde{\BA}_1^F(\theta^0)     \\
     \widetilde{\BA}_1^F(\theta^1)\\
  \end{array}  \right]_D
\right] \nn \\
&&\cdot \left[ \widetilde{\BI}_1 - \mu_0 \left[\begin{array}{c}
    \widetilde{\BI}_0^1(\theta^0)   \\
    \widetilde{\BI}_0^1(\theta^1)\\
  \end{array} \right]\widetilde{\BA}^C_0(2\theta^0)^{-1}\left[\begin{array}{c}
    \widetilde{\BI}_1^0(\theta^0)   \\
    \widetilde{\BI}_1^0(\theta^1)\\
  \end{array} \right]^t\left[ \begin{array}{c}
    \widetilde{\BA}_1^F(\theta^0)\\
     \widetilde{\BA}_1^F(\theta^1)\\
  \end{array}  \right]_D
\right], \label{repre-mc}
\end{eqnarray}
where $\widetilde{\BI}_1$ is $2\times 2$ identity matrix and the
subscript-$_D$ denotes the transformation of a vector into a
diagonal matrix. Similarly, the representation
$\widetilde{\BM}^{SL}_2$ of two-level iteration matrix based on
shifted Laplace operator, can be easily obtained (cf. \cite{CV11}).
For the iteration operator $\BM^{FC}_2$, its representation on
$E^{\theta^0}_{2h_1}$ is given by
\begin{eqnarray}
\widetilde{\BM}^{FC}_2 &=& \left[
\widetilde{\BI}_1 - \mu_1\left(
\widetilde{\BI}_1 -\left[
\begin{array}{c}
    \widetilde{\BS}_1^F(\theta^0)    \\
      \widetilde{\BS}_1^F(\theta^1)\\
  \end{array} \right]_D\right)\right] \nn \\
 && \cdot \left[ \widetilde{\BI}_1 - \mu_0
\left[\begin{array}{c}
    \widetilde{\BI}_0^1(\theta^0)   \\
    \widetilde{\BI}_0^1(\theta^1)\\
  \end{array} \right]\widetilde{\BA}^C_0(2\theta^0)^{-1}\left[\begin{array}{c}
    \widetilde{\BI}_1^0(\theta^0)   \\
    \widetilde{\BI}_1^0(\theta^1)\\
  \end{array} \right]^t\left[ \begin{array}{c}
    \widetilde{\BA}_1^F(\theta^0)\\
     \widetilde{\BA}_1^F(\theta^1)\\
  \end{array}  \right]_D
\right].\label{repre-mfc}
\end{eqnarray}
Then the spectral radius of
$\widetilde{\BM}^{SL}_2,\widetilde{\BM}^{FC}_2$ and
$\widetilde{\BM}^{C}_2$ for different $\theta^0\in \Theta_{\rm low}$
can be obtained analytically and numerically.

\begin{figure}[htbp]
\centering
    \includegraphics[width=2.5in]{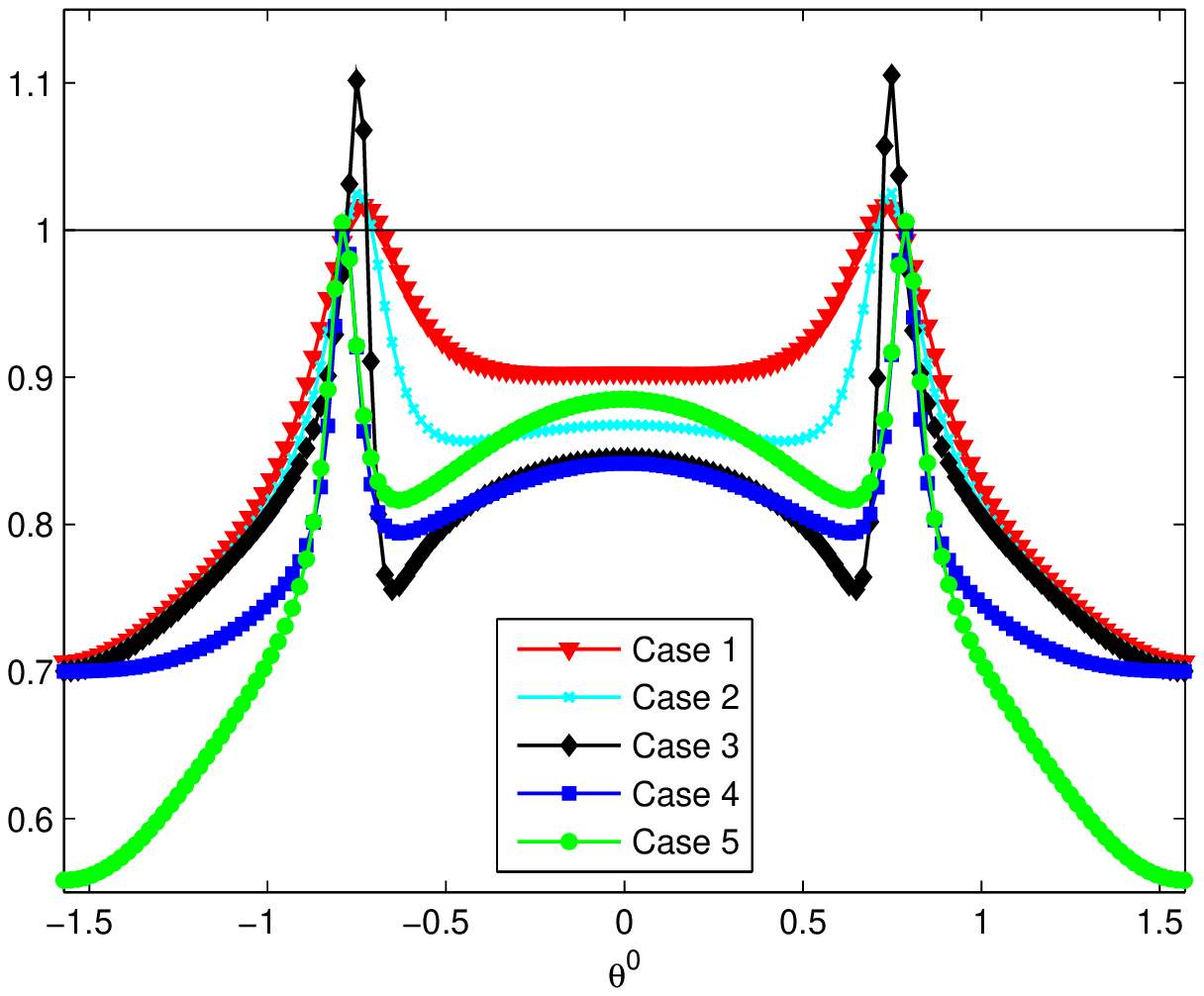}
    \includegraphics[width=2.5in]{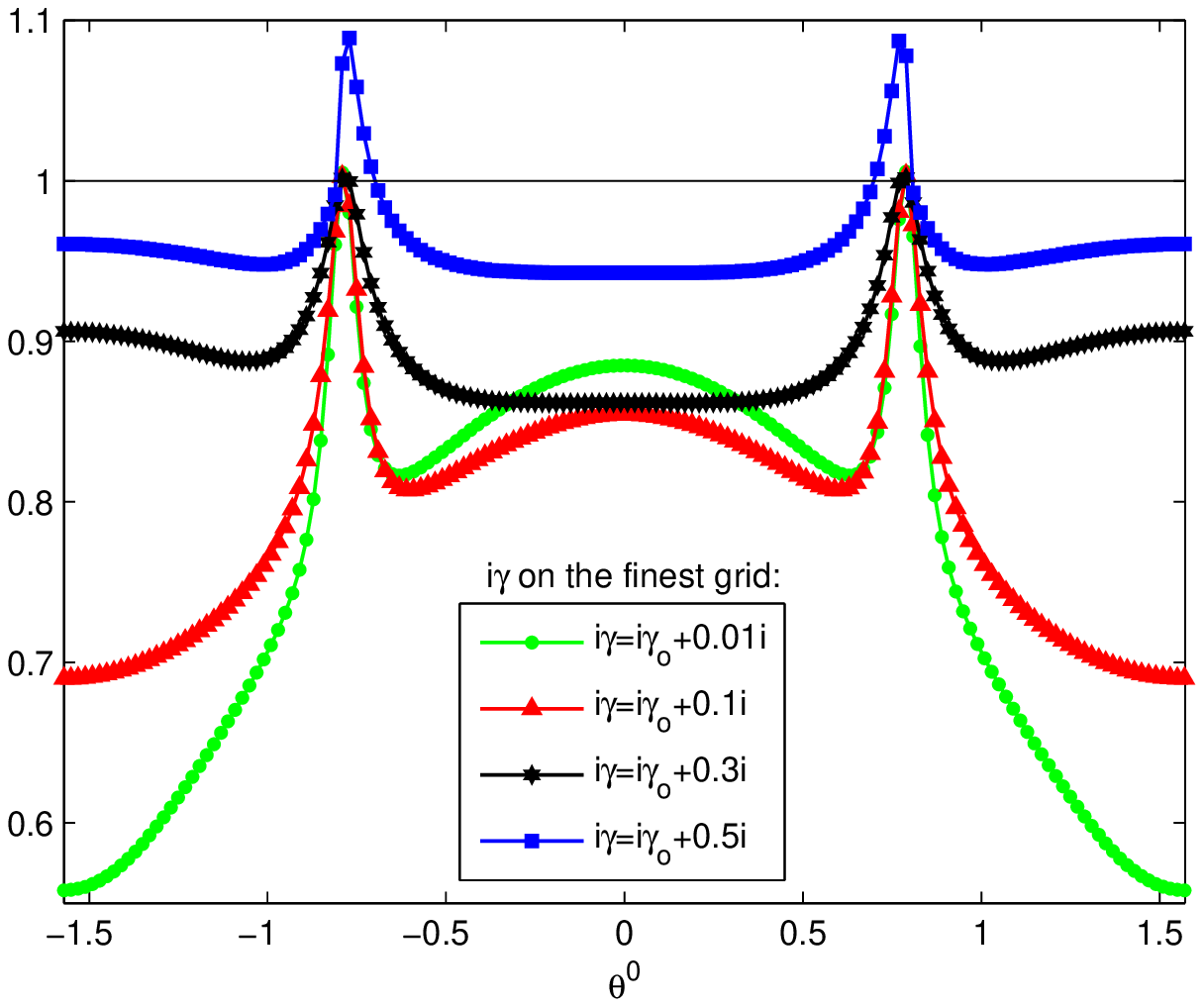}
    \caption{\small Left: Case 1-Case 3: $\rho(\widetilde{\BM}^{SL}_2)$ with $\beta=0.5,0.3 $ and $0.1$ over $\theta^0\in \Theta_{\rm low}$, Case 4
    and Case 5 denote for $\rho(\widetilde{\BM}^{FC}_2)$ and
    $\rho(\widetilde{\BM}^{C}_2)$ with ${\bf i}\gamma = {\bf i}\gamma_o+0.01{\bf
i}$ (${\bf i}\gamma_o\thickapprox -0.085$ when $t=0.8$) on the
finest grid.
    Right: $\rho(\widetilde{\BM}^{C}_2)$ with different ${\bf i}\gamma$. Since $t=kh>1$ on the coarse grid, we
    formally choose the corresponding parameter ${\bf i}\gamma = {\bf i}\gamma_o+0.05{\bf i}$ for the cases in this figure. }\label{fig-jacobi-sl-fc-c}
\end{figure}

The left graph of Figure \ref{fig-jacobi-sl-fc-c} shows the spectral
radius of $\widetilde{\BM}^{SL}_2$ with $\beta=0.5,0.3$ and $0.1$
over $\theta^0\in \Theta_{\rm low}$, and $\widetilde{\BM}^{FC}_2$,
$\widetilde{\BM}^{C}_2$ under $\omega$-JAC relaxation for $t=0.8$ on
the finest grid. In order to make some comparisons between different
algorithms, the penalty parameter ${\bf i}\gamma$ in the CIP-FEM is
always chosen as a complex number in the following if there is no
special notation. If no confusion is possible, we will always denote
$t={\it constant}$ for $t=\kappa h_L$ on the finest grid. The
parameter in $\omega$-JAC relaxation is always chosen as
$\omega=0.6$ in this section. We observe that for most of the
frequencies $\theta^0\in \Theta_{\rm low}$ the spectral radius is
smaller than one. The spectral radius tends to larger than one only
under a few frequencies. Actually, the appearance of such a
resonance is caused by the coarse grid correction and originates
from the inversion of the coarse grid discretization symbol,
$\widetilde{\BA}^C_0(2\theta^0)$ in (\ref{repre-mc}) and
(\ref{repre-mfc}) for instance. The minimum of
$|\widetilde{\BA}^C_0(2\theta^0)|$ maximizes the spectral radius of
$\widetilde{\BM}^{FC}_2$ and $\widetilde{\BM}^{C}_2$ for fixed $t$
and ${\bf i}\gamma$.

Moreover, for different choices of parameter $\beta$ in shifted
Laplace operator, the properties of $\widetilde{\BM}^{FC}_2$ and
$\widetilde{\BM}^{C}_2$ with ${\bf i}\gamma={\bf i}\gamma_o+0.01{\bf
i}$ on the finest grid always perform better than that of
$\widetilde{\BM}^{SL}_2$. Although the spectral radius of
$\widetilde{\BM}^{SL}_2$ with $\beta=0.1$ is almost equivalent to
that of $\widetilde{\BM}^{FC}_2$ over the frequencies spreading near
zero, there is a relatively large resonance causing by the coarse
grid correction. Small $\beta$ in (\ref{shift-lap}) deteriorates the
coarse grid correction. However, large $\beta$ amplifies the
corresponding spectral radius over the frequencies spreading near
zero. We can refer to a recent work \cite{CV11} about the choice of
minimal complex shift parameter in shifted Laplace preconditioner.
The choice of ${\bf i}\gamma$ in $\widetilde{\BM}^{FC}_2$ or
$\widetilde{\BM}^{C}_2$ would also influence its spectral radius.
The right graph of Figure \ref{fig-jacobi-sl-fc-c} shows the
spectral radius of $\widetilde{\BM}^{C}_2$ with different ${\bf
i}\gamma $ under $\omega$-JAC relaxation for $t=0.8$. The influence
of ${\bf i}\gamma$ in $\widetilde{\BM}^{C}_2$ is similar to $\beta$
in shifted Laplace preconditioner. For other small $t<1$, the
properties of $\widetilde{\BM}^{SL}_2$, $\widetilde{\BM}^{FC}_2$ and
$\widetilde{\BM}^{C}_2$ can also be observed as above. Moreover,
when $t=\kappa h$ on the finest grid is small enough, the spectral
radiuses of $\widetilde{\BM}^{SL}_2$, $\widetilde{\BM}^{FC}_2$ and
$\widetilde{\BM}^{C}_2$ can all be smaller than one over all
$\theta^0 \in \Theta_{\rm low}$. Due to the fact that our aim is to
study the influences of smoothing corrections in two- and
three-level methods, in the following we do not always assume $t$ to
be sufficiently small.

\begin{figure}[htbp]
\centering
    \includegraphics[width=2.5in]{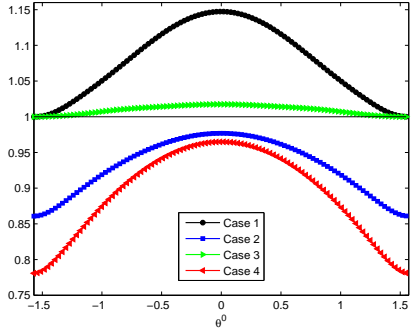}
    \caption{\small Case 1 ($t=\sqrt{3}$) and Case 3 ($t=4$) denote for $\rho(\widetilde{\BM}^{SL}_2)$ with $\beta=0.8$ over $\theta^0 \in \Theta_{\rm low}$, Case 2 ($t=\sqrt{3}$) and Case 4 ($t=4$)
    denote for $\rho(\widetilde{\BM}^{C}_2)$
    with $\gamma = 0.8$.
    }\label{fig-jacobi-sm}
\end{figure}

Alternatively, for large $t$, the maximum of spectral radiuses of
$\widetilde{\BM}^{SL}_2$, $\widetilde{\BM}^{FC}_2$ and
$\widetilde{\BM}^{C}_2$ appear to be spread around $\theta^0=0$.
Actually, for large $t$ the main factor determining the spectral
radius is the smoothing operator rather than the coarse grid
correction. We assume that $\omega$-JAC relaxation is used as
smoother. It can be observed from Figure \ref{fig-jacobi-sm} that
$\rho(\widetilde{\BM}^{SL}_2)$ and $\rho(\widetilde{\BM}^{C}_2)$
reach a maximum in $\theta^0=0$ for $t=\sqrt{3}$ and $t=4$. Indeed,
it has been analyzed in previous section and \cite{CV11} that
$\rho(\BS^J_h)$ for CIP-FEM and shifted Laplace preconditioner
(\ref{shift-lap}) with continuous piecewise P1 approximation reaches
a maximum at $\theta=0$ or $\theta=\pi$. From (\ref{sfem-sm}), it is
easy to see that the $\omega$-JAC relaxation for discrete system
from standard P1 FEM approximation is inefficient when $t$ is near
$\sqrt{3}$. Thus, two-level method $\BM^{FC}_2$ is inefficient in
this case. However, ${\BM}^{SL}_2$ and ${\BM}^{C}_2$ can still be
applied due to the adding of complex parts on the fine and coarse
grids. One can observe from Figure \ref{fig-jacobi-sm} that the
spectral radiuses of $\rho(\widetilde{\BM}^{SL}_2)$ and
$\rho(\widetilde{\BM}^{C}_2)$ is always larger than one when
$t=\sqrt{3}$, which implies the divergence of the associated
two-level methods when applying $\omega$-JAC smoother. However, the
spectral radius is smaller than one when $t=4$, which may permit the
use of $\omega$-JAC as a smoother on very coarse grids, but we shall
not use this and utilize GMRES smoothing on coarse grids instead.

Furthermore, to get a comprehensive insight into the influence of
multiple coarse grid corrections, we next carry out a three-level
local Fourier analysis. The iteration operator of three-level method
by nonsymmetric version of Algorithm \ref{modi-mu} is derived as
$(I-T_2)(I-T_1)(I-T_0)$. Similar to the two-level method, when
CIP-FEM is used for smoothing on $\Ct_2,\Ct_1$ and $\Ct_0$, the
iteration matrix can be deduced to be
\[
\BM^C_3 = \Big(\BI_2 - \mu_2 (\BI_2 - \BS^C_2)(\BA^C_2)^{-1} \BA^F_2\Big)
\Big(\BI_2 - \mu_1 \BI^2_1 (\BI_1 - \BS^C_1)(\BA^C_1)^{-1}\BI^1_2
\BA^F_2\Big)\Big(\BI_2 - \mu_0 \BI^2_0(\BA^C_0)^{-1} \BI^0_2
\BA^F_2\Big).
\]
Define the four dimensional $4h$-harmonics by
\[
E_{4h_2}^{\theta^0} :={\rm span} \{ \varphi_{h_2}(\theta^{00},x) ,
\varphi_{h_2}(\theta^{01},x) , \varphi_{h_2}(\theta^{10},x),
\varphi_{h_2}(\theta^{11},x) \},
\]
where $\theta^{\alpha0} = \frac{\theta^{\alpha}}{2},\theta^{\alpha1}
=
 \frac{\theta^{\alpha}}{2}-{\rm sign}( \frac{\theta^{\alpha}}{2})\pi,\alpha=0,1$. For any
$\theta^0 \in \Theta_{\rm low}$, the three-level operators leaves
the space of $4h$-harmonics $E_{4h_2}^{\theta^0} $ invariant (cf.
\cite{WJ04}). This yields a block diagonal representation of
$\BM^C_3$ with the following $4\times 4$ matrix
$\widetilde{\BM}^C_3$:

\begin{eqnarray}
\widetilde{\BM}^C_3 &=& \left[ \widetilde{\BI}_2 - \mu_2\left(
\widetilde{\BI}_2 -\widetilde{\BS}^C_2\right)
(\widetilde{\BA}^C_2)^{-1} \widetilde{\BA}^F_2
\right] \nn \\
&&\cdot \left[ \widetilde{\BI}_2 - \mu_1 \widetilde{\BI}_1^2
 \left(
\widetilde{\BI}_1 -\left[
\begin{array}{c}
    \widetilde{\BS}_1^C(\theta^0)    \\
      \widetilde{\BS}_1^C(\theta^1)\\
  \end{array} \right]_D\right) \left[  \begin{array}{cc}
    \widetilde{\BA}_1^C(\theta^0)   \\
     \widetilde{\BA}_1^C(\theta^1)\\
  \end{array}  \right]_D^{-1} (\widetilde{\BI}_2^1)^t\widetilde{\BA}^F_2
\right]\label{repre-mc-3} \\
&& \cdot \left[ \widetilde{\BI}_2 - \mu_0 \widetilde{\BI}_1^2
\widetilde{\BI}_0^1 \widetilde{\BA}^C_0(2\theta^0)^{-1}
(\widetilde{\BI}_1^0)^t (\widetilde{\BI}_2^1)^t \widetilde{\BA}^F_2
\right],\nn
\end{eqnarray}
where $\widetilde{\BI}_2 $ is $4\times 4$ identity matrix,
$\widetilde{\BS}^C_2  = \left[
\begin{array}{cc}
    \widetilde{\BS}_2^C(\theta^{00})   \\
    \widetilde{\BS}_2^C(\theta^{01})    \\
      \widetilde{\BS}_2^C(\theta^{10})    \\
     \widetilde{\BS}_2^C(\theta^{11})    \\
  \end{array} \right]_D$, $\widetilde{\BA}^C_2
= \left[  \begin{array}{cc}
    \widetilde{\BA}_2^C(\theta^{00})   \\
     \widetilde{\BA}_2^C(\theta^{01})\\
     \widetilde{\BA}_2^C(\theta^{10})\\
     \widetilde{\BA}_2^C(\theta^{11})\\
  \end{array}  \right]_D$,
$\widetilde{\BA}^F_2  = \left[  \begin{array}{cc}
    \widetilde{\BA}_2^F(\theta^{00})     \\
    \widetilde{\BA}_2^F(\theta^{01})     \\
    \widetilde{\BA}_2^F(\theta^{10})     \\
    \widetilde{\BA}_2^F(\theta^{11})     \\
  \end{array}  \right]_D$,
$\widetilde{\BI}_1^2  = \left[
\begin{array}{cc}
    \widetilde{\BI}_1^2(\theta^{00})  \quad 0 \\
    \widetilde{\BI}_1^2(\theta^{01})  \quad 0  \\
     0 \quad \widetilde{\BI}_1^2(\theta^{10})    \\
     0 \quad \widetilde{\BI}_1^2(\theta^{11})    \\
  \end{array} \right]$, $\widetilde{\BI}_0^1 = \left[\begin{array}{c}
    \widetilde{\BI}_1^0(\theta^0)   \\
    \widetilde{\BI}_1^0(\theta^1)\\
  \end{array} \right]$, and $\widetilde{\BI}_2^1 ,\widetilde{\BI}_1^0  $ are defined similarly.
For another approach, we assume that standard FEM is applied on
$\Ct_2$ for smoothing and CIP-FEM is used for correction on $\Ct_1$
and $\Ct_0$, then the associated iteration matrix is given by
\[
\BM^{FC}_3 =  \Big(\BI_2- \mu_2(\BI_2- \BS^F_2)   \Big) \Big(\BI_2 - \mu_1 \BI^2_1 (\BI_1 -
\BS^C_1)(\BA^C_1)^{-1}\BI^1_2 \BA^F_2\Big)\Big(\BI_2 - \mu_0
\BI^2_0(\BA^C_0)^{-1} \BI^0_2 \BA^F_2\Big).
\]
Moreover, similar to $\BM^{SL}_2$, we denote by $\BM^{SL}_3$ the
three-level iteration matrix, which takes the smoothing based on
shifted Laplace operator (\ref{shift-lap}) with standard FEM
approximation. Then the representations $\widetilde{\BM}^{FC}_3$,
$\widetilde{\BM}^{SL}_3$ with respect to $\BM^{FC}_3$ and
$\BM^{SL}_3$ respectively on the $4h$-harmonics
$E_{4h_2}^{\theta^0}$ can be derived directly.

\begin{figure}[htbp]
\centering
    \includegraphics[width=2.5in]{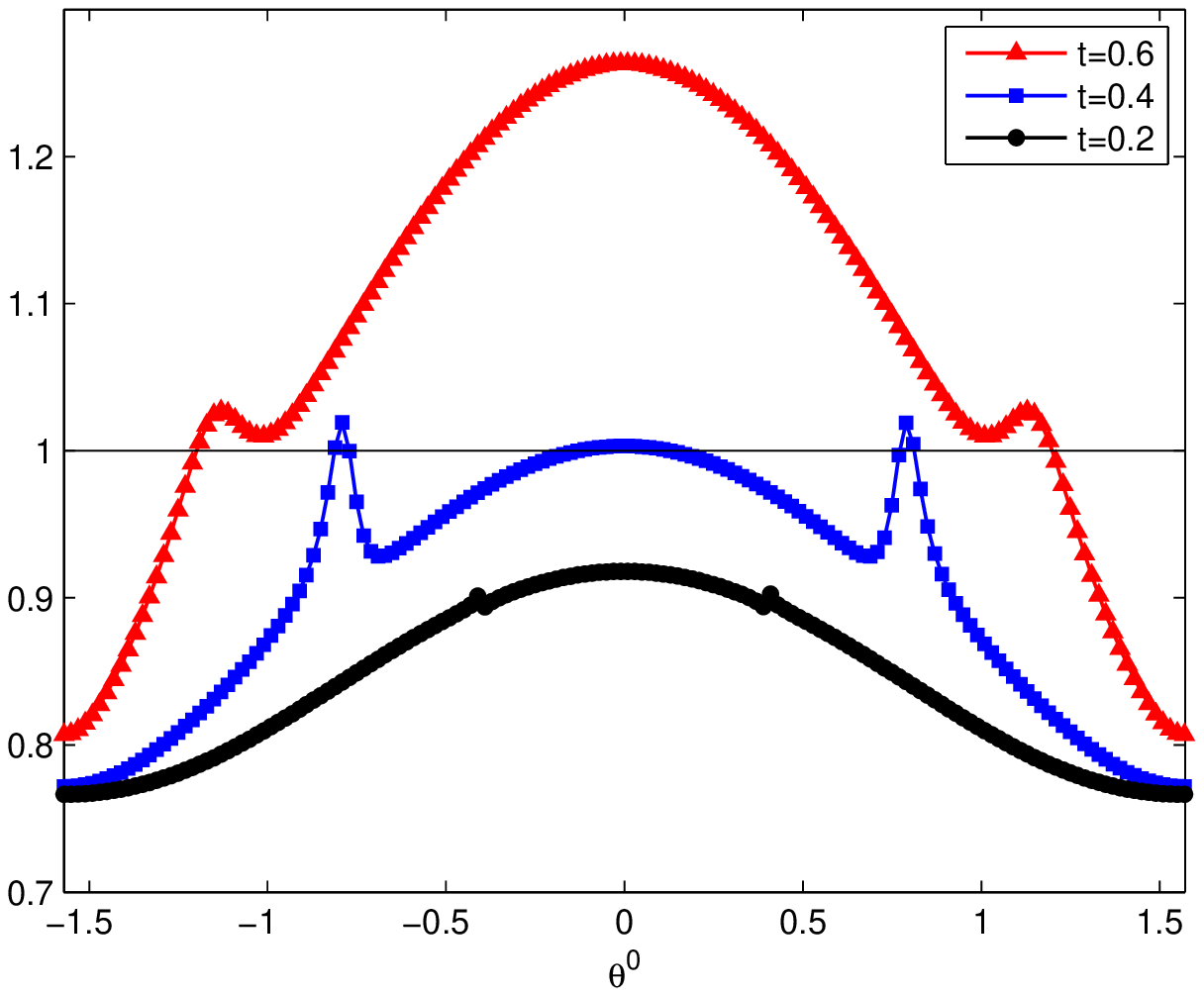}
    \includegraphics[width=2.5in]{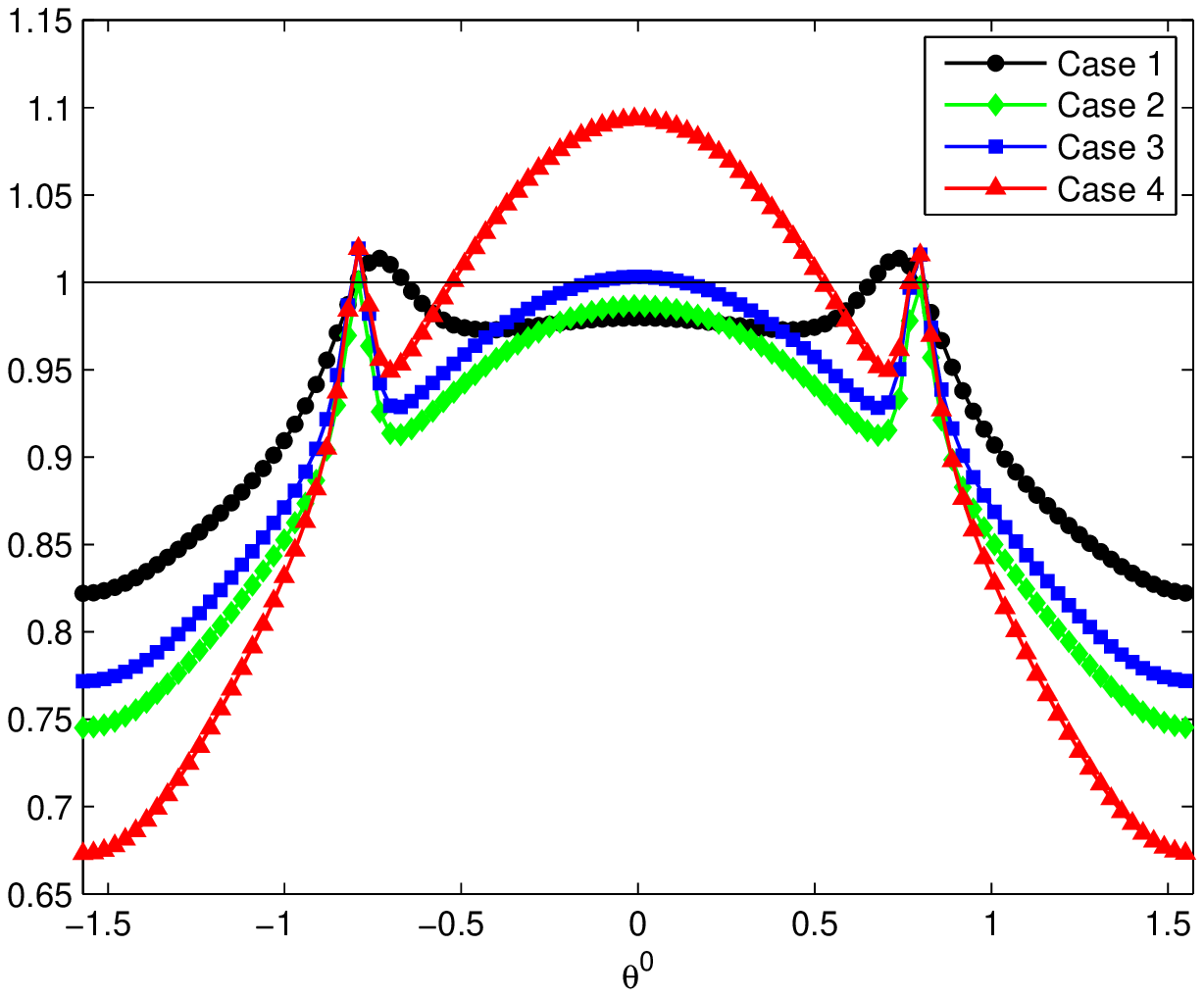}
    \caption{\small Left: $\rho(\widetilde{\BM}^C_3)$ over $\theta^0 \in \Theta_{\rm low}$ for the cases $t=0.6,0.4,0.2$.
    Right ($t=0.4$): Case 1-Case 2: $\rho(\widetilde{\BM}^{SL}_3)$,
    $\rho(\widetilde{\BM}^{C}_3)$, Case 3-Case 4:
    $\rho(\widetilde{\BM}^{FC}_3)$ with one step and two step intermediate coarse grid (the 2nd level) correction.
    $\omega$-JAC relaxation is applied for smoothing. The parameters in shifted Laplace operator and CIP-FEM are chosen as $\beta=0.5$, ${\bf i}
    \gamma = {\bf i}\gamma_o + 0.01{\bf i}$ when $t<1$ and ${\bf i}
    \gamma = {\bf i}\gamma_o + 0.05{\bf i}$ when $t>1$.
    }\label{fig-jacobi-3d}
\end{figure}

For fixed wave number $\kappa$, the performance of three-level
methods is shown in Figure \ref{fig-jacobi-3d}. For instance, one
can observe from the left graph of Figure \ref{fig-jacobi-3d} that a
coarser initial grid used for $\BM^C_3$ will deteriorate its
convergence. The right graph of Figure \ref{fig-jacobi-3d} shows the
spectral radiuses for three kinds of approaches: $\BM^{SL}_3$,
$\BM^{C}_3$ and $\BM^{FC}_3$. For $t=0.4$, both of $\BM^{FC}_3$ and
$\BM^{C}_3$ perform better than $\BM^{SL}_3$ over most of the
frequencies, and the spectral radius of $\BM^{FC}_3$ is similar to
that of $\BM^{C}_3$ over most of frequencies in this case. Actually,
the performance of $\BM^{FC}_3$ is always comparable with
$\BM^{C}_3$. This is also true for two-level method, which can be
observed from the left graph of Figure \ref{fig-jacobi-sl-fc-c}.
Thus, in order to reduce the computational cost, the modified
multilevel methods (Algorithm \ref{modi-mu}) with stable CIP-FEM
corrections on coarser grids and standard FEM corrections on finer
grids can be usually utilized in practical. From the right graph of
Figure \ref{fig-jacobi-3d}, we can also see that the convergence of
$\BM^{FC}_3$ does not always perform better with more steps of
intermediate coarse grid correction.

\section{Numerical results}\label{sec-numer}

In this section, we present some numerical results to illustrate the
performance of Algorithm \ref{alg1} and modified multilevel method
(Algorithm \ref{modi-mu}) for two Helmholtz problems in two
dimension. The multilevel method is used as a preconditioner in
outer GMRES iterations (PGMRES). Gauss-Seidel relaxation is always
used as smoother when $l \in S_L$, and the smoothing steps are
always chosen as one step if there is no any annotation. At the
$l$-th level, the discrete problem reads ${\bf A}_l {\bf u}_l = {\bf
F}_l$. We denote by $ {\bf r}^n_l = {\bf F}_l - {\bf A}_l {\bf
u}^n_l$ the residual with respect to the $n$-th iteration. The
PGMRES algorithm terminates when
\begin{equation}\label{e:6}
\|{\bf r}^n_l\|/\|{\bf r}^0_l\| \leq 10^{-6} .\nn
\end{equation}
The number of iteration steps required to achieve the desired accuracy is denoted by {\bf iter}.

\begin{exm}
{\rm We consider a two dimensional Helmholtz equation with the first order
absorbing boundary condition (cf. \cite{Wu09,Wu12}):
\begin{align*}
-\triangle u - \kappa^2 u &= f:= \frac{\sin(\kappa r)}{r}\quad {\rm in} \
\Omega,  \\
\frac{\partial u}{\partial n} + {\ii} \kappa u &= g \quad {\rm on} \
\partial \Omega.
\end{align*}
Here $\Omega$ is a unit square with center $(0,0)$ and $g$ is chosen
such that the exact solution is
\[
u = \frac{\cos(\kappa r)}{\kappa} - \frac{\cos \kappa+ \ii \sin
\kappa}{\kappa(J_0(\kappa)+{\ii} J_1(\kappa))} J_0 (\kappa r),
\]
where $J_{\nu}(z)$ are Bessel functions of the first kind.

\smallskip

\begin{table}[htbp]
\caption{\footnotesize Iteration counts of PGMRES for discrete
system from CIP-P1 and CIP-P2 ($\kappa=100$). The smoothing
relaxations on fine and coarse grids are all based on CIP-FEM.}
\begin{center}
\footnotesize
\begin{tabular}{|c|cccc|}
\hline
Level & 2 & 3 & 4 & 5 \\
\hline
DOFs & 16641 & 66049 & 263169 & 1050625 \\
\hline
iter (P1)  & 27 & 24 & 21 & 20\\
\hline
iter (P2)  & 26& 22 & 19 & 18\\
\hline
\end{tabular}\label{ex1-1}
\end{center}
\end{table}

In this example, we assume that the coarsest level of multilevel
method satisfies $\kappa h_0/p \thickapprox 2$ for $\kappa\leq 360$.
For $400 \leq \kappa \leq 600$, we choose the coarsest grid
condition such that $\kappa h_0/p\approx 1.1 \thicksim 1.7$. The
parameters in (\ref{Jump-var}) are chosen as $\gamma_e\equiv
0.01+0.07{\bf i}$ (cf. \cite{Wu09}) for CIP-P1, and $\gamma_e\equiv
0.005+0.035{\bf i}$ for piecewise P2 CIP-FEM (CIP-P2). We first test
the algorithms for the case with wave number $\kappa=100$. When
Algorithm \ref{alg1} is applied to the discrete system
$A^C_Lu_L=F_L$, the smoothing relaxations are all based on the
CIP-FEM approximation on fine and coarse grids. Table \ref{ex1-1}
shows the corresponding iteration counts of PGMRES for discrete
system from CIP-P1 and CIP-P2 approximations. We can observe that
for fixed $\kappa$ the algorithm is robust on different levels.

\begin{table}[htbp]
\caption{\footnotesize Iteration counts of PGMRES for discrete
system from FEM-P1. The smoothing relaxations on fine and coarse
grids are all based on FEM-P1 for $\kappa=100,200$. }\label{ex1-2d}
\begin{center}
\footnotesize
\begin{tabular}{|c|c|c|cccc|}
\hline
\multirow{4}{*}{$\kappa=100$} & \multicolumn{2}{|c|}{Level}   &    2    &   3    &  4     &  5  \\
\cline{2-7}
 & \multicolumn{2}{|c|}{DOFs}                               & 16641 & 66049 & 263169 & 1050625 \\
\cline{2-7}
 & \multicolumn{2}{|c|}{iter ($m_2=1$)}                               & 58 & 87 & 94 & 90\\
 \cline{2-7}
 & \multicolumn{2}{|c|}{iter ($m_2=20$)}                                & 46& 53 & 49 & 47\\
\hline \hline

\multirow{4}{*}{$\kappa=200$} & \multicolumn{2}{|c|}{Level}    &    2    &   3    &  4     &  5  \\
\cline{2-7}
 & \multicolumn{2}{|c|}{DOFs}                                 & 66049 & 263169 & 1050625 & 4198401 \\
\cline{2-7}
 & \multicolumn{2}{|c|}{iter ($m_2$=1)}                       & 177   &  $>200$   & $>200$ & $>200$   \\
\hline

\end{tabular}
\end{center}
\end{table}

For fixed $\kappa$, when the grid is fine enough to get the
accuracy, standard FEM can be utilized again to discretize the
problem. But standard multigrid method fails to solve the
corresponding discrete system for the case with large wave number.
In the following, we mainly apply the multilevel method presented in
Algorithm \ref{modi-mu} with different smoothing strategies. Table
\ref{ex1-2d} shows the iteration counts of this
multilevel-preconditioned GMRES method with GMRES smoothing on
coarse grids for $\kappa=100,200$, and the smoothing relaxations on
fine and coarse grids are all based on standard P1 FEM (FEM-P1)
approximation. We find that the iteration count is mesh independent
for fixed $\kappa$, but it increases rapidly with larger wave
number. For instance, for the case $\kappa=200$ with $m_2=1$, the
iteration counts of PGMRES will be more than 200, even when the
GMRES smoothing is performed by $m_2=20$ steps, the iteration counts
are still more than 100. Although more steps of GMRES smoothing can
reduce the total iteration counts, it requires more memory to store
data in the computation. Actually, the slow convergence of the
algorithm mainly lies in the bad approximation on coarse grids.
Next, we apply CIP-FEM to construct stable coarse grid corrections.

\begin{figure}[htbp]
\centering
    \includegraphics[width=2.5in]{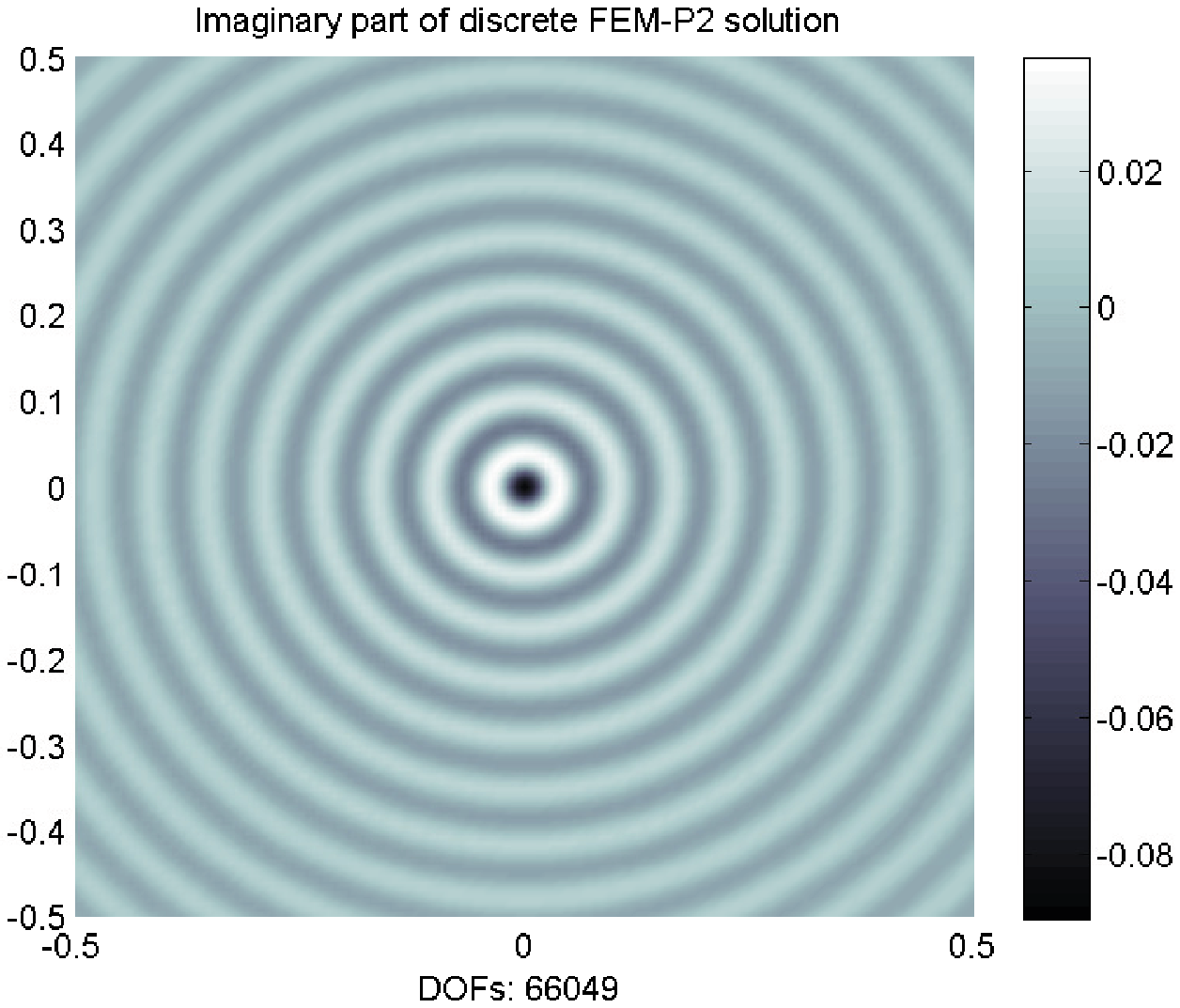}
    \includegraphics[width=2.5in]{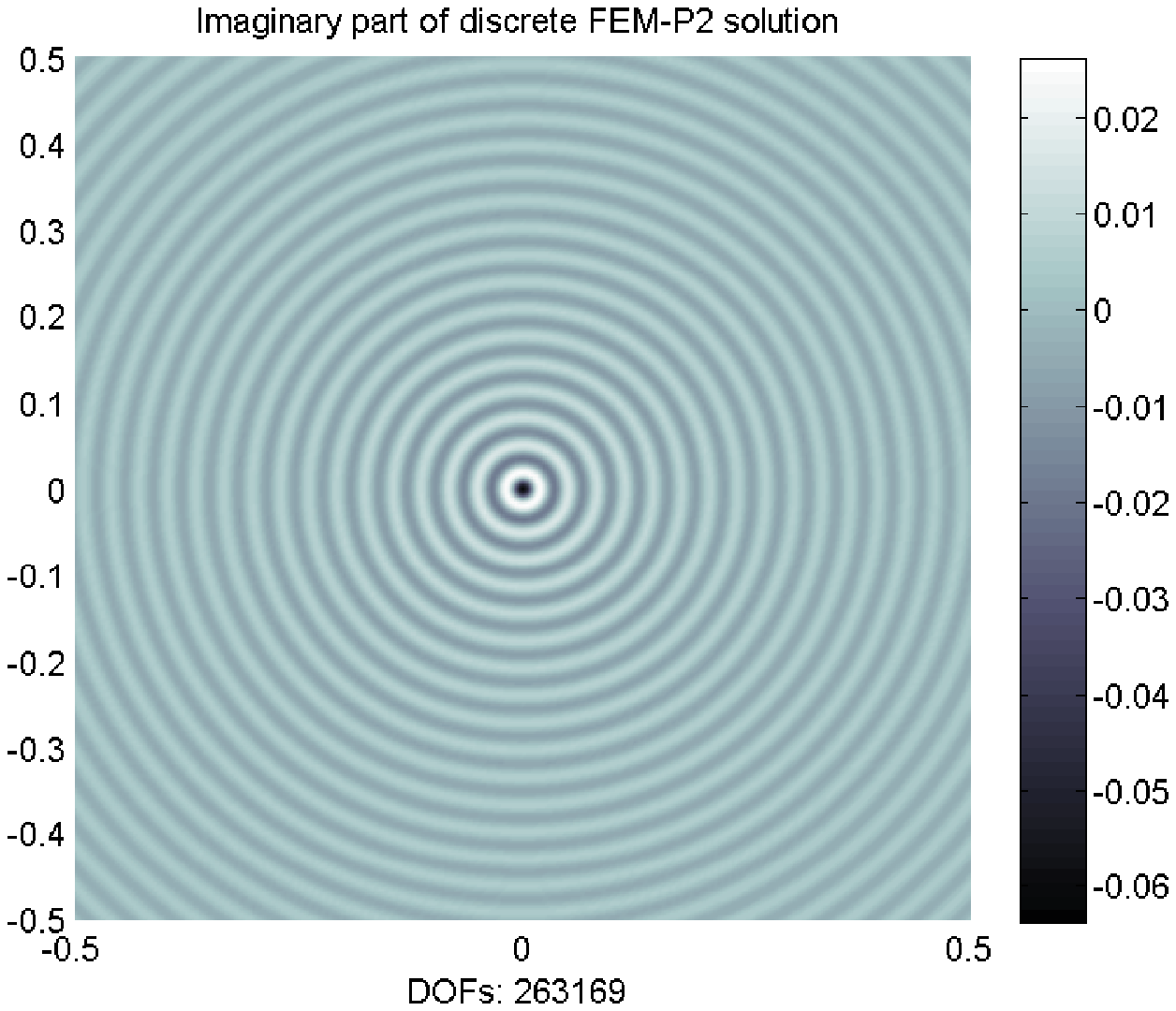}
    \caption{\footnotesize Surface plot of imaginary part of discrete FEM-P2 solutions for $\kappa=100$ (left) and $\kappa=200$ (right)
    on the grid with mesh condition $\kappa h / p \approx 0.55$.}\label{ex1-fig-plot}
\end{figure}

\begin{table}[htbp]\footnotesize
\caption{\footnotesize Iteration counts of PGMRES for discrete
system from FEM-P1 approximation ($\kappa=100$). Case 1: the
smoothing relaxations on fine and coarse grids are both based on
shifted Laplace operator (\ref{shift-lap}) with $\beta=0.2$  by
FEM-P1, Case 2: the smoothing relaxations on fine and coarse grids
are both based on CIP-P1.
}\label{ex1-3d}
\begin{center}
\begin{tabular}{|c|cccc|}
\hline
Level & 2 & 3 & 4 & 5 \\
\hline
DOFs & 16641 & 66049 & 263169 & 1050625 \\
\hline
iter (Case 1)  & 52 & 64 & 66 & 63\\
\hline
iter (Case 2)  & 27& 24 &  22 & 21 \\
\hline
\end{tabular}
\end{center}
\end{table}

Figure \ref{ex1-fig-plot} displays the surface plot of imaginary
part of discrete FEM-P2 solutions for $\kappa=100,200$ on the grid
with mesh condition $\kappa h / p \approx 0.55$. Indeed, the
discrete solutions have correct shapes and amplitudes as the exact
solutions. Table \ref{ex1-3d} presents the comparing of shifted
Laplace operator with FEM-P1 and the original Helmholtz operator
with CIP-P1. For the case $\kappa=100$, we can see that the second
approach has the minimum iteration counts. The iteration counts
shown in Table \ref{ex1-4} examine the performance of PGMRES based
on Algorithm \ref{modi-mu} with different steps of GMRES smoothing.
The CIP-P1 and CIP-P2 are only used for approximations on coarse
grids. From Tables \ref{ex1-1}, \ref{ex1-3d} and \ref{ex1-4}, it
suggests that Algorithm \ref{modi-mu} is also very efficient as
Algorithm \ref{alg1} and the second approach in Table \ref{ex1-3d}
which uses CIP-FEM on both fine and coarse grids, and the growth in
GMRES smoothing steps does not always reduce iteration counts
obviously. This supports the similar phenomena shown in the right
graph of Figure \ref{fig-jacobi-3d}. Hence, in the following we will
only use one step ($m_2=1$) of GMRES smoothing in Algorithm
\ref{modi-mu}.

\begin{table}[htbp]\footnotesize
\caption{\footnotesize Iteration counts of PGMRES based on Algorithm
\ref{modi-mu} for discrete system from FEM-P1 and FEM-P2 with
different steps of post GMRES smoothing
($\kappa=100$).}\label{ex1-4}
\begin{center}
\begin{tabular}{|c|ccc|}
\hline
Level  & 3 & 4 & 5 \\
\hline
DOFs   & 66049 & 263169 & 1050625 \\
\hline
iter (P1, $m_2=1$)   & 24 & 23 & 22\\
\hline
iter (P1, $m_2=10$)   & 21  & 23  & 23 \\
\hline
iter (P2, $m_2=1$)    & 27  & 21   & 19  \\
\hline
iter (P2, $m_2=10$)     & 32 &   21 & 18  \\
\hline
\end{tabular}
\end{center}
\end{table}

\begin{table}[htbp]
\caption{\footnotesize Iteration counts of PGMRES based on Algorithm
\ref{modi-mu} for discrete system from FEM-P1 and FEM-P2 for the
cases $\kappa=50,200,360$ with the coarsest grid condition $\kappa
h_0 / p \thickapprox 2$.}\label{ex1-5}
\begin{center}
\footnotesize
\begin{tabular}{|c|c|c|ccc|}
\hline
\multirow{4}{*}{$\kappa=50$} & \multicolumn{2}{|c|}{Level}      &   3    &  4     &  5  \\
\cline{2-6}
 & \multicolumn{2}{|c|}{DOFs}                                  &  16641 & 66049 & 263169   \\
\cline{2-6}
 & \multicolumn{2}{|c|}{iter (P1)}                                   &   15     &  15     &   15      \\
 \cline{2-6}
 & \multicolumn{2}{|c|}{iter (P2)}                                    &   23   &   14  &     13   \\
\hline \hline

\multirow{4}{*}{$\kappa=200$} & \multicolumn{2}{|c|}{Level}       &   3    &  4     &  5  \\
\cline{2-6}
 & \multicolumn{2}{|c|}{DOFs}                                  & 263169 & 1050625 & 4198401 \\
\cline{2-6}
 & \multicolumn{2}{|c|}{iter (P1)}                               &  49  &  57 &  55  \\
\cline{2-6}
 & \multicolumn{2}{|c|}{iter (P2)}                                &  44 &  41 & 36 \\
\hline \hline

\multirow{4}{*}{$\kappa=360$} & \multicolumn{2}{|c|}{Level}   &    2    &   3    &  4      \\
\cline{2-6}
 & \multicolumn{2}{|c|}{DOFs}                                 &    263169     &  1050625      &   4198401      \\
\cline{2-6}
 & \multicolumn{2}{|c|}{iter (P1)}                     &  111    &   115   &  112    \\
\cline{2-6}
 & \multicolumn{2}{|c|}{iter (P2)}                     &   41   & 44     &  39    \\
\hline


\end{tabular}
\end{center}
\end{table}

\begin{table}[htbp]
\caption{\footnotesize Iteration counts of PGMRES based on Algorithm
\ref{modi-mu} for discrete system from FEM-P1 and FEM-P2 for the
cases $\kappa=400,500,600$. The coarsest grid condition is chosen
with mesh size $h_0 \thickapprox 0.00276$ such that $\kappa
h_0/p\thickapprox 1.1 \thicksim 1.7$.}\label{ex1-6}
\begin{center}
\footnotesize
\begin{tabular}{|c|c|c|cc | c|c|c|cc|}
\hline
\multirow{5}{*}{FEM-P1} & \multicolumn{2}{|c|}{Level}   &    2    &   3     &  \multirow{5}{*}{FEM-P2} & \multicolumn{2}{|c|}{Level}    &    2    &   3        \\
\cline{2-5} \cline{7-10}
 & \multicolumn{2}{|c|}{DOFs}                                 & 1050625 & 4198401  &  & \multicolumn{2}{|c|}{DOFs}                                  & 1050625 & 4198401 \\
\cline{2-5} \cline{7-10}
 & \multicolumn{2}{|c|}{iter ($\kappa=400$)}                 &    30   &   26     &      & \multicolumn{2}{|c|}{iter ($\kappa=400$)}                  & 21  &   14        \\
 \cline{2-5} \cline{7-10}
 & \multicolumn{2}{|c|}{iter ($\kappa=500$)}                 &   71   &   50     &       & \multicolumn{2}{|c|}{iter ($\kappa=500$)}                &  29  &  25      \\
  \cline{2-5} \cline{7-10}
 & \multicolumn{2}{|c|}{iter ($\kappa=600$)}                         &  156    &  141  &       & \multicolumn{2}{|c|}{iter ($\kappa=600$)}                 &   43   &     47      \\
\hline
%

\end{tabular}
\end{center}
\end{table}

Tables \ref{ex1-5} and \ref{ex1-6} show the iteration counts of
PGMRES based on Algorithm \ref{modi-mu} with CIP-FEM used for
correction only on coarse grids for discrete system from FEM-P1 and
FEM-P2. For fixed $\kappa$, one can observe that the iteration
counts are robust on different levels. Due to the same coarsest grid
used for $\kappa = 400,500,600$, the growth of iteration counts is a
little faster than linear when $\kappa \geq 400$. But when the
FEM-P2 is applied in particular, the growth of iteration counts with
increasing wave number is more stable than FEM-P1.

}
\end{exm}

\bigskip

\begin{exm}
{\rm In a unit square domain $\Omega$ with center $(0,0)$,
we consider the Helmholtz problem (\ref{Hem})-(\ref{Robin}) with discontinuous wave number, which is defined by
\[
\kappa =
\begin{cases}
\kappa_1, & \text{if}\;\; (x_1,x_2)\in  (-0.5,0)\times (0,0.5) \bigcup (0,0.5)\times (-0.5,0) ,\\
\kappa_2, & \text{elsewhere},
\end{cases}
\]
where $\kappa_2 = q \kappa_1, q>1$. We set the Robin boundary condition (\ref{Robin}) with $g=0$ and
the external force $f(x)$ to be a narrow Gaussian point source (cf. \cite{Engquist2}) located at
$(r_1,r_2) = (-0.25,-0.25)$:
\[
f(x_1,x_2) = e^{ -(\frac{4\kappa}{\pi})^2 ( (x_1-r_1)^2 + (x_2-r_2)^2 )  }.
\]

\begin{figure}[h]
\centering
    \includegraphics[width=2.5in,height=2.1in]{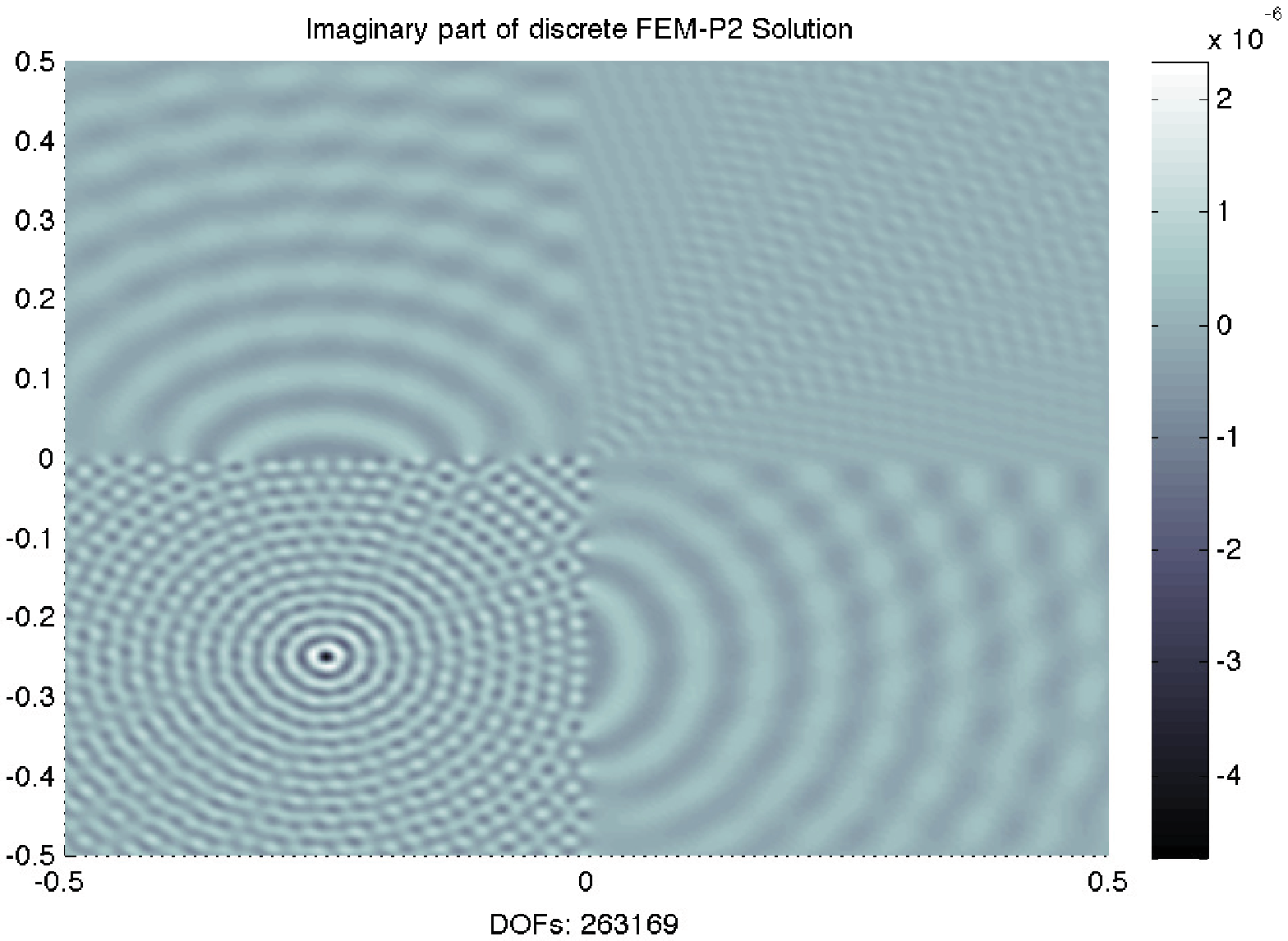}
    \includegraphics[width=2.5in,height=2.1in]{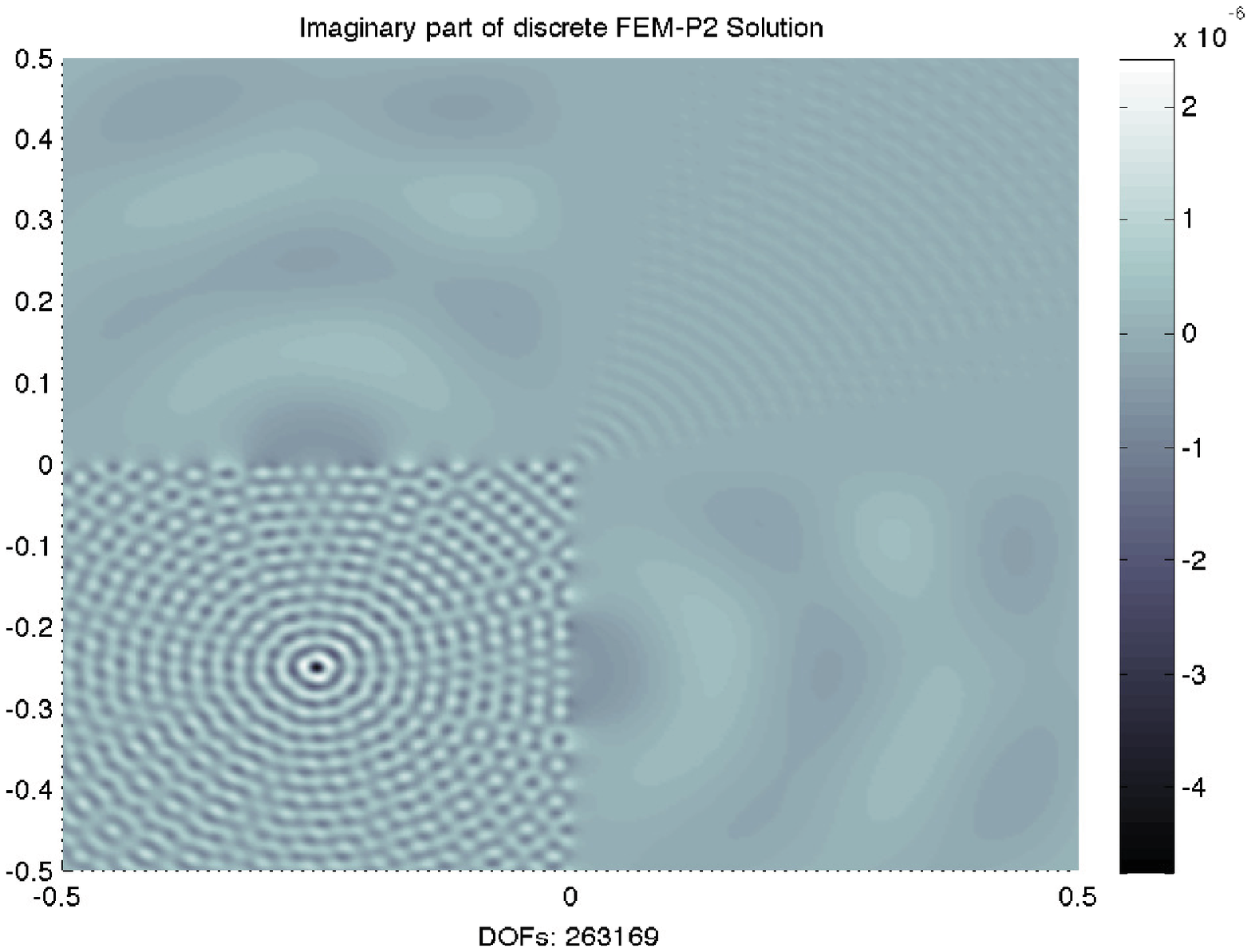}
    \caption{\footnotesize Surface plot of imaginary part of discrete FEM-P2 solutions for $\kappa_2=300,q=3$ (left) and $\kappa_2=300,q=10$ (right)
    on the grid with mesh condition $\kappa_2 h / p \approx 0.8$. }\label{ex2-fig}
\end{figure}

\begin{table}[h]
\caption{\footnotesize Iteration counts of PGMRES based on Algorithm
\ref{modi-mu} for discrete system from FEM-P1 and FEM-P2 for the
cases $\kappa_2=180 \ (q=3,10),\ \kappa_2=300 \
(q=3,10)$.}\label{ex2-1}
\begin{center}
\footnotesize
\begin{tabular}{|c|c|c|cc|}
\hline
\multirow{6}{*}{$\kappa_2=180$} & \multicolumn{2}{|c|}{Level}      &   3   & 4      \\
\cline{2-5}
 & \multicolumn{2}{|c|}{DOFs}                            &        263169     &  1050625    \\
\cline{2-5}
 & \multicolumn{2}{|c|}{iter (P1,$q=3$)}                   &    26    &      27      \\
 \cline{2-5}
 & \multicolumn{2}{|c|}{iter (P1,$q=10$)}               &    28    &  29   \\
  \cline{2-5}
 & \multicolumn{2}{|c|}{iter (P2,$q=3$)}                  &   16     &  15     \\
  \cline{2-5}
 & \multicolumn{2}{|c|}{iter (P2,$q=10$)}                &     17   &    15  \\
\hline \hline

\multirow{6}{*}{$\kappa_2=300$} & \multicolumn{2}{|c|}{Level}     &   3   & 4      \\
\cline{2-5}
 & \multicolumn{2}{|c|}{DOFs}                                    &  1050625      &   4198401 \\
\cline{2-5}
 & \multicolumn{2}{|c|}{iter (P1,$q=3$)}                 &      30   &     30     \\
 \cline{2-5}
 & \multicolumn{2}{|c|}{iter (P1,$q=10$)}                 &    32     &    33  \\
  \cline{2-5}
 & \multicolumn{2}{|c|}{iter (P2,$q=3$)}                  &   16   &    15  \\
  \cline{2-5}
 & \multicolumn{2}{|c|}{iter (P2,$q=10$)}                 &   16     &  15   \\
\hline

\end{tabular}
\end{center}
\end{table}

In this example, we test Algorithm \ref{modi-mu} for the Helmholtz
problem with discontinuous wave number. Due to the fact that
$\kappa_2 = \max\{\kappa_1,\kappa_2\}$, the coarsest mesh condition
is according to the choice of $\kappa_2 h_0 / p$. Figure
\ref{ex2-fig} displays the surface plot of imaginary part of
discrete FEM-P2 solutions for $\kappa_2=300$ with $q=3$ and $q=10$
on the grid with mesh condition $\kappa_2 h / p \approx 0.8$. The
iteration counts of PGMRES based on Algorithm \ref{modi-mu} for
different $\kappa_2$ and $q$ are listed in Table \ref{ex2-1}. For
fixed $\kappa_2$ and polynomial order, we note that the iteration
counts are robust on different levels. Similar to the first example,
the PGMRES iteration for FEM-P2 is more stable than FEM-P1.

}
\end{exm}


\begin{thebibliography}{bib}
\small

\bibitem{Adams}
R. Adams, Sobolev Spaces, Academic Press, New York, 1975.

\bibitem{Brandt77}
A. Brandt, Multi-level adaptive solutions to boundary-value
problems, Math. Comp., 31 (1977), pp. 333--390.

\bibitem{Brandt97}
A. Brandt and I. Livshits, Wave-ray multigrid method for
standing wave equations, Electron. Trans. Numer. Anal., 6 (1997),
pp. 162--181.

\bibitem{CV11}
S. Cools and W. Vanroose, Local Fourier Analysis of the Complex
Shifted Laplacian preconditioner for Helmholtz problems, preprint,
2011.

\bibitem{EEO01}
H.C. Elman, O.G. Ernst, and D.P. O'Leary, A multigrid method
enhanced by Krylov subspace iteration for discrete Helmholtz
equations, SIAM J. Sci. Comput., 23 (2001), pp. 1291--1315.


\bibitem{Engquist1}
B. Engquist and L. Ying, Sweeping preconditioner for the
Helmholtz equation: hierarchical matrix representation, Comm. Pure
Appl. Math., 64 (2011), pp. 697--735.

\bibitem{Engquist2}
B. Engquist and L. Ying, Sweeping preconditioner for the
Helmholtz equation: moving perfectly matched layers, Multiscale Model. Simul.,
9 (2011), pp. 686--710.


\bibitem{Erlangga04}
Y. A. Erlangga, C. Vuik, and C. W. Oosterlee, On a class of
preconditioners for solving the Helmholtz equation, Appl. Numer.
Math., 50 (2004), pp. 409--425.

\bibitem{Erlangga06}
Y. A. Erlangga, C. W. Oosterlee, and C. Vuik, A novel multigrid
based preconditioner for heterogeneous Helmholtz problems, SIAM J.
Sci. Comput., 27 (2006), pp. 1471--1492.



\bibitem{Erlangga08}
Y.A. Erlangga, Advances in iterative methods and preconditioners for
the Helmholtz equation, Arch. Comput. Methods Eng., 15 (2008), pp.
37--66.

\bibitem{Ernst11}
O.G. Ernst and M.J. Gander, Why it is difficult to solve Helmholtz
problems with classical iterative methods?, in: I. Graham, T. Hou,
O. Lakkis, R. Scheichl (Eds.), Numerical Analysis of Multiscale
Problems, Springer, 2011.




\bibitem{Wu09}
X. Feng and H. Wu, Discontinuous Galerkin methods for the Helmholtz
equation with large wave number, SIAM J. Numer. Anal, {\bf 47}
(2009), pp. 2872--2896.

\bibitem{Wu09hp}
X. Feng and H. Wu, hp-discontinuous Galerkin methods for the
Helmholtz equation with large wave number, Math. Comp, {\bf 80}
(2011), pp. 1997--2024.

\bibitem{Ihlenburg_book}
F. Ihlenburg, Finite Element Analysis of Acoustic Scattering,
Springer-Verlag, New York, 1998.


\bibitem{Kim}
S. Kim and S. Kim, Multigrid simulations for high-frequency
solutions of the Helmholtz problem in heterogeneous media, SIAM J.
Sci. Comput., 24 (2002), pp. 684--701.

\bibitem{Livshits}
I. Livshits and A. Brandt, Accuracy properties of the wave-ray
multigrid algorithm for Helmholtz equations, SIAM J. Sci. Comput.,
28 (2006), pp. 1228--1251.

\bibitem{VW}
P.S. Vassilevski and J. Wang, An application of the abstract
multilevel theory to nonconforming finite element methods, {\it SIAM
J. Numer. Anal.}, {\bf 32} (1995), 235--248.

\bibitem{WJ04}
R. Wienands and W. Joppich, Practical Fourier Analysis for Multigrid
Methods, Chapman \& Hall/CRC, London, 2004.

\bibitem{Wu12}
H. Wu, Pre-asymptotic error analysis of CIP-FEM and FEM for
Helmholtz equation with high wave number. Part I: Linear version, to
appear, 2013.

\bibitem{Wu12-1d}
L. Zhu, E. Burman, and H. Wu, Continuous Interior Penalty Finite
Element Method for Helmholtz Equation with Large Wave Number: One
Dimension Analysis, to appear, 2013.

\bibitem{Wu12-hp}
L. Zhu and H. Wu, Pre-asymptotic error analysis of CIP-FEM and FEM
for Helmholtz equation with high wave number. Part {II}: $hp$
version, SIAM J. Numer. Anal., to appear, 2013. (See also arXiv:
http://arxiv.org/pdf /1204.5061v1).


\end{thebibliography}
\end{document}